\documentclass[11pt]{amsart}
\usepackage{generalwriting}

\usepackage{tikz-cd}
\microtypesetup{tracking=true}
\UseMicrotypeSet[tracking]{alltext-nott}
\SetTracking{encoding=*, family=ppl}{10}

\usepackage[dvipsnames]{xcolor}

\begin{document}
\thanks{linxuan.li@qmul.ac.uk}
\title[Tropical Regularity]{Tropical Regularity}
\author{Linxuan Li}

\begin{abstract}
We introduce the notion of tropical regularity as a generalization of tropical smoothness and present a geometric classification conjecture for tropical regular cycles. Roughly, the conjecture predicts that every regular tuple either forms a regular fibration or descends via a lattice projection to one of finitely many minimal-model cycles. Our main results establish this conjecture in three fundamental cases:
unmixed regular tuples $(T,\ldots,T,[F])$ with $L(T)=\{0\}$ and $[F]$ spanning its ambient space; regular tropicalizations of complete intersections with generic coefficients; and regular tropical $2$-cycles satisfying the Hodge index inequality
and spanning their ambient space.
\end{abstract}

\maketitle
\tableofcontents

\section{Introduction}

Over the past 30 years, tropical geometry has developed rapidly
at the interface of algebraic geometry and combinatorial geometry.
One of its main contributions is the translation of difficult problems
in algebraic geometry into combinatorial problems.
Examples include Mikhalkin's enumeration of curves on toric surfaces
through weighted counts of tropical curves \cite{Mikhalkin05}
and the tropical proof of the Brill--Noether theorem through divisor
theory on metric graphs \cite{CoolsDraismaPayneRobeva12}.

At the same time, tropical geometry has developed its own world.
These include tropical intersection theory \cite{AllermannRau10}, cohomologically tropical varieties \cite{AksnesAminiPiquerezShaw25},
and the theory of tropical ideals \cite{MaclaganRincon18}.
These theories are formulated directly in terms of tropical objects
and apply beyond those admitting algebraic realizations.
Tropical geometry has thereby become a subject in its own right,
with its own structures, methods, and questions.

Both aspects are reflected in the notion of a \emph{tropical manifold},
which adapts the classical idea of a manifold to the tropical setting
using combinatorial local models
\cite{MikhalkinZharkov14,MikhalkinRau}.
Informally, tropical manifolds can be viewed as tropical spaces
without points of multiplicity greater than one.
The key property of tropical manifolds is \emph{tropical smoothness}.
In particular, a tropical fan-cycle is tropical smooth if and only if
it is $\mathbb Z$-linearly isomorphic to a Bergman cycle
\cite{MikhalkinRau}.

Fink's characterization of Bergman fans gives a explicit description
of this condition:
\begin{theorem}[\cite{AlexFink}]\label{AlexFink}
A tropical $k$-fan $F$ in $\mathbb R^n$ is a Bergman $k$-fan if and only if
$$
\deg(\underbrace{\max(0,x_1,\dots,x_n)\cdots \max(0,x_1,\dots,x_n)}_{k\text{ times}}\cdot [F])=1
$$
where $x_i$'s form the standard basis of $(\R^n)^\vee$.
\end{theorem}

Thus tropical smoothness of a fan-cycle is equivalent to this
numerical condition after a suitable change of lattice coordinates. Moreover, this numerical condition suggests a broader notion,
which we call \emph{tropical regularity}.
A $k$-dimensional tropical fan-cycle $[F]$ in $\mathbb R^n$
is called \emph{regular} if there exist $k$ tropical regular functions
$T_1,\ldots,T_k$, not necessarily distinct, such that
\[
    \deg\bigl(T_1\cdots T_k\cdot[F]\bigr)=1.
\]
This condition relaxes Fink's criterion by allowing arbitrary tropical
regular functions in place of repeated copies of the standard function $H$.
In particular, every tropical smooth fan-cycle is regular.
This leads to the central question of this paper.

\begin{question}\label{Q}
    \textit{How far between tropical regularity and tropical smoothness?}
\end{question}

In intersection theory,
Adiprasito, Huh, and Katz proved that the real Chow ring
$A^*(M)_{\mathbb R}$ of any loopless matroid satisfies
Poincar\'e duality, the hard Lefschetz theorem
and the Hodge--Riemann relations \cite{AHK18}.
Their theorem requires no realizability assumption and proves
the log-concavity of the absolute values of the coefficients
of the characteristic polynomial of $M$. The matroid Chow ring is canonically identified with
$A^*(X_{\Sigma_M})$, where $X_{\Sigma_M}$ is the toric variety
associated with the fine Bergman fan $\Sigma_M$
\cite{FeichtnerYuzvinsky04,AHK18}.
It is therefore natural to ask which tropical regular cycles admit suitable refinements whose Chow rings satisfy analogous duality
and positivity properties. This opens a direction in the intersection theory of tropical regular cycles that lies beyond the scope of this paper.

In combinatorial geometry, tropical regular cycles form
a much larger class than tropical smooth cycles.
The difference is already visible for one-dimensional fans
in $\mathbb R^2$, as illustrated in the picture below:

\begin{figure}[htbp]
    \centering
    \begin{tikzpicture}[>=stealth, line width=0.8pt, font=\small]
        \begin{scope}
            \draw[->] (0,0) -- (1.45,0)
                node[right] {$(1,0)$};
            \draw[->] (0,0) -- (-1.15,1.15)
                node[above left] {$(-1,1)$};
            \draw[->] (0,0) -- (0,-1.3)
                node[below] {$(0,-1)$};
            \fill (0,0) circle (1.5pt);
            \node[align=center] at (0,-2.2)
                {Tropical smooth $1$-fan $C$};
        \end{scope}

        \begin{scope}[xshift=6cm]
            \draw[->] (0,0) -- (1.45,0)
                node[right] {$(1,0)$};
            \draw[->] (0,0) -- (0,1.65)
                node[above] {$(0,1)$};
            \draw[->] (0,0) -- (-0.75,1.5)
                node[above left] {$(-1,q)$};
            \draw[->] (0,0) -- (0,-1.3)
                node[below] {$(0,-1)$};
            \node[right] at (0,-0.65) {$q+1$};
            \fill (0,0) circle (1.5pt);
            \node[align=center] at (0,-2.2)
                {Regular $1$-fan $F_q$, $q\in\mathbb Z_{>1}$\\
                 not tropical smooth};
        \end{scope}
    \end{tikzpicture}
    \label{fig:regular-versus-smooth}
\end{figure}

In this paper, we focus on comparing the combinatorial geometry
of tropical regular cycles with that of tropical smooth cycles.
In dimension one, a complete classification was obtained
in \cite{Li25}.
It constructs a canonical lattice projection from every regular
$1$-cycle to a direct sum of Bergman $1$-fans, along which every
tropical regular function of intersection number one descends
\cite[Theorem~3.16]{Li25}.
This gives a positive answer to Question~\ref{Q} in dimension one.

Building on this geometric classification, Alexander Esterov raised the following conjectural classification of higher-dimensional tropical regular cycles satisfying the Hodge index inequality (cf. Definition \ref{Hodgeindex}).

\begin{conjecture*}[A.Esterov]\label{trop}
     There exists a finite class $\mathcal{MD}_n$ of {\it minimal model} $k$-cycles, such that every regular $k$-cycle $[F]$ in $\R^n$ satisfying the Hodge index inequality and $\mathrm{span}_{\R}|F|=\R^n$, is one of the following: for any regular tuple $(T_1,\dots,T_k,[F])$,

\begin{itemize}
    \item[\textbf{I}:] either the regular tuple $(T_1,\dots,T_k,[F])$ forms a regular fibration (cf. Definition \ref{fibration}),
    
    \item[\textbf{II}:] or there exists a saturated $\Z$-linear projection $\pi\colon\R^n\twoheadrightarrow\R^m$ such that $\pi_*[F]\in\mathcal{MD}_n$. Moreover, there exists a regular tuple $(M_1,\dots,M_k,\pi_*[F])$ such that
    $$T_i=\pi^*M_i$$
    for all $1\le i\le k$.\\
\end{itemize}
\end{conjecture*}

\begin{remark} 
The Hodge index inequality is satisfied by all of the most important classes of tropical fans, such as matroid fans (see \cite{AHK18}) and tropicalizations of algebraic varieties (see \cite{EsterovECI}).

Fans failing to satisfy the Hodge index inequality do not look ubiquitous in algebraic or combinatorial geometry, or even easy to construct: for instance, the key achievement of the fundamental paper \cite{BabaeeHuh19} was to construct at least one such ``irreducible'' fan in $\R^4$, as it gave a counterexample to a certain strong version of the Hodge conjecture.
\end{remark}

\begin{remark}
Conjecture~\hyperref[trop]{A.~Esterov} does not directly answer
Question~\ref{Q}.
For instance, we prove the conjecture holds for tropical regular
$2$-cycles spanning $\mathbb R^n$ and satisfying the Hodge index inequality.
However, the relationship between the resulting minimal-model class
$\mathcal{MD}_n$ and tropical smooth cycles remains to be understood.
This comparison is an interesting problem in its own right.
\end{remark}

\begin{remark}
When the stable intersection is a point with weight $0$, it is closely related to Esterov's question on reducible engineered complete intersections \cite{EsterovECI}. In that setting, one is led to study pairs of tropical fans $(F_1,F_2)$ with $F_1\cdot F_2=0$, and Esterov observes that this becomes an incidence problem for the projectivized spans of their maximal cones. Our Theorem \ref{TF=0_thm} gives the corresponding unmixed divisor-theoretic form: $\deg(T^k\cdot [F])=0$ precisely when the lineality space $L(T)$ meets $\langle\sigma\rangle$ nontrivially for every maximal cone $\sigma$ of $F$.
\end{remark}

\begin{remark}
The classifications of degree-zero and degree-one intersections in this paper provide useful evidence for geometric questions on tropical complete intersections. In particular, they are closely related to Question 5.6 in \cite{EsterovSCI}. In the setting considered there, a tropical complete intersection may be written as $(\Sigma_{P_1},\ldots,\Sigma_{P_k},[F])$, where $[F]$ is a tropical $k$-cycle in $\mathbb{R}^n$, and its Euler characteristic is defined by the formal expression
$$\frac{\Sigma_{P_k}}{1+\Sigma_{P_k}}\cdots\frac{\Sigma_{P_1}}{1+\Sigma_{P_1}}\cdot [F],$$
where each quotient $\Sigma_{P_i}/(1+\Sigma_{P_i})$ is interpreted via its Taylor expansion. Question 5.6 asks, among other things, for a classification of tropical complete intersections with Euler characteristic of small absolute value, including the cases $0$ and $1$. From this perspective, our classification of the cases $\deg(T^k\cdot [F])=0$ and $\deg(T^k\cdot [F])=1$ gives an unmixed intersection-theoretic model for the first two extremal cases of Esterov's question.
\end{remark}

\subsection{Structure and results}

We now describe the structure of the paper and state its main results.
Section~\ref{S2} reviews the necessary background in tropical geometry, especially the tropical intersection theory.
Section~\ref{S3} introduces the notion of tropical regularity and
formulates Conjecture~\hyperref[trop]{A.~Esterov}.

In Section~\ref{S4}, we prove the conjecture for tropical regular
$2$-cycles satisfying the Hodge index inequality and spanning their
ambient space.

\begin{theorem}[cf.Theorem \ref{Thm_1}]
There exists a finite collection $\mathcal{MD}_n$ of
\emph{minimal model} $2$-cycles with the following property.
Let $[F]$ be a regular $2$-cycle in $\R^n$ satisfying the Hodge index
inequality and $\operatorname{span}_{\R}|F|=\R^n$.
For every regular tuple $(T_1,T_2,[F])$, one of the following holds:
\begin{itemize}
    \item[\textbf{I}:]
    The tuple $(T_1,T_2,[F])$ forms a regular fibration.

    \item[\textbf{II}:]
    There exists a saturated $\Z$-linear projection
    $\pi\colon\R^n\twoheadrightarrow\R^m$ such that
    $\pi_*[F]\in\mathcal{MD}_n$.
    Moreover, there is a regular tuple $(M_1,M_2,\pi_*[F])$ satisfying
    \[
        T_i=\pi^*M_i
        \qquad (i=1,2).
    \]
\end{itemize}
\end{theorem}

In Section~\ref{S5}, we prove the conjecture for regular
tropicalizations of complete intersections with generic coefficients.
In this case, the resulting models or fibration bases are Bergman cycles.

\begin{theorem}[cf.Theorem \ref{Thm_2}]
Let
\[
    U=V(f_1,\ldots,f_{n-k})\subseteq(\C^*)^n
\]
be a $k$-dimensional very affine variety with generic coefficients.
If $[F]:=[\Trop(U)]$ is regular, then every regular tuple
$(T_1,\ldots,T_k,[F])$ determines a saturated $\Z$-linear projection
$\pi\colon\R^n\twoheadrightarrow\R^m$ such that one of the following holds:
\begin{enumerate}
    \item
    The tuple $(T_1,\ldots,T_k,[F])$ forms a regular fibration along
    $\pi$, and the set-theoretic image $\pi(|F|)$ is the support
    of a Bergman fan in $\R^m$.

    \item
    The pushforward $\pi_*[F]$ is a Bergman $k$-cycle in $\R^m$,
    and there exist tropical regular functions
    $M_1,\ldots,M_k\in\mathrm{TR}(\R^m)$ satisfying
    \[
        T_i=\pi^*M_i
        \qquad (1\le i\le k),
        \qquad
        \deg\bigl(M_1\cdots M_k\cdot\pi_*[F]\bigr)=1.
    \]
\end{enumerate}
\end{theorem}

Finally, Section~\ref{S6} studies unmixed intersection numbers
$\deg(T^k\cdot[F])$, focusing on the values zero and one.
We first give a geometric criterion for the case of degree $0$.

\begin{theorem}[cf.Theorem \ref{TF=0_thm}]
Let $F$ be a tropical $k$-fan in $\R^n$ and let $T$ be a tropical
regular function.
Then $\deg(T^k\cdot[F])=0$ if and only if
\[
    L(T)\cap\langle\sigma\rangle\neq\{0\}
    \qquad\text{for every }\sigma\in F^{(k)},
\]
where $\langle\sigma\rangle:=\operatorname{span}_{\R}(\sigma)$.
\end{theorem}

For the case of degree $1$, we obtain the following rigidity theorem,
extending Fink's characterization to an arbitrary tropical regular
function with trivial lineality space.

\begin{theorem}[cf.Theorem \ref{Thm_3}]
Let $[F]$ be a tropical $k$-cycle in $\R^n$ satisfying
$\operatorname{span}_{\R}|F|=\R^n$.
Suppose that $T$ is a tropical regular function such that
\[
    \deg(T^k\cdot[F])=1,
    \qquad
    L(T)=\{0\},
\]
where $L(T)$ is its lineality space
(cf.~Definition~\ref{TR}).
Then there exists a lattice automorphism
$\psi\colon\R^n\to\R^n$ such that:
\begin{enumerate}
    \item
    $\psi_*[F]$ is a Bergman $k$-cycle spanning $\R^n$;

    \item
    there exists a tropical regular function
    $M\in\mathrm{TR}(\R^n)$ satisfying
    \[
        T=\psi^*M,
        \qquad
        \deg\bigl(M^k\cdot\psi_*[F]\bigr)=1.
    \]
\end{enumerate}
Moreover, $\mathrm{Newt}(T)$ is an $n$-simplex.
\end{theorem}

\paragraph*{\textbf{Acknowledgements.}}
The author is deeply grateful to \textsc{Alexander Esterov} and \textsc{Alex Fink} for their patient supervision, fruitful discussions, and valuable suggestions.  \\

\section{Tropical Intersection Theory}\label{S2}
Tropical intersection products of tropical cycles are central objects in this paper. In this section, we recall the basic preliminaries on tropical intersection theory that will be used throughout.

\begin{definition}[Tropical fans]\label{tropical_fans}
A {tropical fan} of dimension $k$ in $\mathbb{R}^n$ is a pair $(F,\omega_F)$ in which $F$ is a rational polyhedral fan in $\mathbb{R}^n$ pure of dimension $k$ and $\omega_F\colon F^{(k)}\to\mathbb{Z}_{>0}$ is a weight function satisfying the {balancing condition}: for every $\tau\in F^{(k-1)}$,
$$\sum_{\sigma\in F^{(k)},\,\tau<\sigma}\omega_F(\sigma)\,v_{\sigma/\tau}\in\langle\tau\rangle,$$
where $v_{\sigma/\tau}\in\mathbb{Z}_\sigma\colon=\mathrm{span}(\sigma)\cap\mathbb{Z}^n$ is the primitive integral generator of $\mathbb{Z}_\sigma/\mathbb{Z}_\tau$ pointing into $\sigma$ (cf. \cite[Section 2]{MikhalkinRau}). When the weights are clear from context we will simply write $F$ for $(F,\omega_F)$.
\end{definition}

Readers unfamiliar with cones and polyhedral fans may consult \cite{Cox} for further details. We begin with a special class of tropical fans called \textit{(tropical) hypersurfaces}:

\begin{example}[hypersurfaces]\label{ex:normal-fan}
Let $P\subset\mathbb{R}^n$ be a full-dimensional lattice polytope, and let $\mathcal{N}(P)$ be its outer normal fan, defined by $C_F:=\mathrm{cl}\{w\in\mathbb{R}^n:\mathrm{face}_w(P)=F\}$ for each face $F$ of $P$ (cf. \cite[Chapter 2]{MikhalkinRau}, \cite{MaclaganSturmfels15}, \cite{Cox}). The duality $\dim C_F+\dim F=n$ identifies the codimension-one cones of $\mathcal{N}(P)$ with the edges of $P$. Let ${\Sigma_P}$ denote the union of all $C_F$ with $\dim F\ge 1$, and weight each maximal cone $C_e$ by the lattice length $\ell(e):=\#(e\cap\mathbb{Z}^n)-1$. Then $\bigl({\Sigma_P},\omega)$ is a tropical fan of pure dimension $n-1$ called \textit{hypersurface}, and the balancing condition at a codimension-two cone $C_F$ holds automatically; cf. \cite{MikhalkinRau} or \cite[Lemma~3.4.6]{MaclaganSturmfels15} for the details.
\end{example}

We now introduce the class of functions used in this paper which arise from lattice polytopes:

\begin{definition}[Tropical regular functions]\label{TR}
A {tropical regular function} on $\mathbb{R}^n$ is a function $T\colon\mathbb{R}^n\to\mathbb{R}$ of the form
$$T(x) \;=\; \max\bigl(\ell_1(x),\dots,\ell_k(x)\bigr), \qquad \ell_i\in(\mathbb{Z}^n)^\vee.$$
Moreover, the lineality space $L(T)$ of $T$ is defined by the set $\{x\in \R^n\mid \ell_i(x)=\ell_j(x)\ \forall 1\le i<j\le k\}$.  We write a tropical regular function simply in the form $\max(\ell_1,\dots,\ell_k)$ for convenience.
\end{definition}

\begin{remark}\label{P_T}
Conversely, every tropical regular function
$T:=\max(\ell_1(x),\dots,\ell_k(x))$
determines a Newton polytope
$\mathrm{Newt}(T):=\mathrm{conv}\{\ell_1,\dots,\ell_k\}.$
As an immediate consequence of Example~\ref{ex:normal-fan}, the polytope
$\mathrm{Newt}(T)$ determines a tropical hypersurface 
$\Sigma_{\mathrm{Newt}(T)}$ denoted by $V(T)$.
\end{remark}

Let $(X,\omega_X)$ and $(Y,\omega_Y)$ be two tropical fans. We will use the word `refinement' in the sense of \cite{AllermannRau10}; note that this differs from the way it is often used in polyhedral geometry. We call $(X,\omega_X)$ a refinement of $(Y,\omega_Y)$ if they have the same support, every cone of $X$ is a subset of a cone of $Y$ and the weight function $\omega_X$ is induced from $\omega_Y$. Two tropical fans are equivalent if they have a common refinement. Readers may consult \cite{AllermannRau10} for more details. 

\begin{definition}[Tropical cycles]
Let $(X,\omega_X)$ be a tropical fan of dimension $k$. We denote the equivalence class under the equivalence relation above by $[(X,\omega_X)]$. We refer to $Z^{\mathrm{aff}}_{k}(\mathbb{R}^n)$ as the set of equivalence classes of tropical fans over $\mathbb{R}^n$. The elements of $Z^{\mathrm{aff}}_k(\mathbb{R}^n)$ are called (affine tropical) k-cycles. Again, we write $[X]$ instead of $[(X,\omega_X)]$ if the weight function is clear.
\end{definition}

\begin{definition}[Tropical Intersection Product]\label{TCI}
    Let $T$ be a tropical regular function and let $[(F,\omega_F)]$ be a $k$-cycle with a representative $(F,\omega_F)$ such that $T$ is linear on each cone $\sigma\in F$.  The tropical intersection product is a $(k-1)$-cycle $T\cdot [F]:=[(\bigcup_{i=0}^{k-1}F^{(i)},\omega_{T\cdot[F]})]\in Z^{\mathrm{aff}}_{k-1}(\mathbb{R}^{n})$:
    \begin{align*}
        \omega_{T\cdot [F]}: F^{(k-1)} &\longrightarrow \mathbb{Z},\\
    \tau &\mapsto \sum_{\sigma\in F^{(k)},\tau<\sigma}T_\sigma(\omega_F(\sigma)v_{\sigma/\tau})-T_\tau(\sum_{\sigma\in F^{(k)},\tau<\sigma}\omega_F(\sigma)v_{\sigma/\tau})
    \end{align*}
    where $T_\sigma:\langle\sigma\rangle\to\mathbb{R}$ is the $\mathbb{Z}$-linear functional defined by $T_\sigma(x)=T(x)$ for all $x\in \sigma$ (similar for $T_\tau$).
\end{definition}

Section 3 of \cite{AllermannRau10} shows that the definition is well-defined: it is independent of the choice of the normal vector $v_{\sigma/\tau}$ and of the representative $F$ of the cycle $[F]$.

\begin{proposition}[\cite{AllermannRau10}]\label{propTCI}
    The tropical intersection product is bilinear and commutative.
    More precisely, let $[C]$ and $[D]$ be $k$-cycles in
    $\mathbb R^n$, with $k\geq 2$, and let $T_1,T_2$ be
    tropical regular functions on $\mathbb R^n$. Then:
    \begin{enumerate}
        \item
        $(T_1+T_2)\cdot[C]
        =T_1\cdot[C]+T_2\cdot[C]$;

        \item
        $T_1\cdot([C]+[D])
        =T_1\cdot[C]+T_1\cdot[D]$;

        \item
        $T_1\cdot T_2\cdot[C]
        =T_2\cdot T_1\cdot[C]$;

        \item
        $\zeta\cdot(T_1\cdot [C])=(\zeta\cdot T_1)\cdot[C]$ for any $\zeta\in\Z_{>0}$.
    \end{enumerate}
    Moreover, $\ell\cdot[C]=0$ for every integral linear form
    $\ell\in(\mathbb Z^n)^\vee$.
\end{proposition}

\begin{remark}\label{TR(Rn)}
    Every tropical regular function
    $T=\max(\ell_0,\ell_1,\dots,\ell_r)$ can be normalized to
    \[
        \widetilde T:=T-\ell_0
        =\max(0,\ell_1-\ell_0,\dots,\ell_r-\ell_0).
    \]
    The function $\widetilde T$ is non-negative and, by Proposition~\ref{propTCI}, has the same intersection product with every tropical cycle as $T$.
    We therefore assume throughout this paper that tropical regular functions are normalized in this way, and write them as $T=\max(0,\ell_1,\dots,\ell_r)$. The family of them is denoted by $\mathrm{TR}(\R^n)$ if $\ell_i\in (\Z^n)^\vee$.
\end{remark}

We will also use the standard intersection product of tropical fans. 
For two tropical fans $X,Y \subset \mathbb{R}^n$ of complementary or expected dimensions, their intersection product $$
[X \cap_{\mathrm{st}} Y],
$$ is the weighted tropical cycle obtained by the fan-displacement rule: set-theoretically, one intersects $X$ with a sufficiently small generic translate of $Y$, and the weights are determined by the usual lattice-index multiplicities. We refer to Fulton--Sturmfels \cite{FultonSturmfels97} for the precise construction in terms of Minkowski weights, and to \cite{MaclaganSturmfels15} for its formulation in tropical intersection theory.

When $[F]\in Z_k^{\mathrm{aff}}(\mathbb{R}^n)$ and $T_1,T_2,\dots,T_k$ are tropical regular functions on $\mathbb{R}^n$, the iterated product
$$
T_1\cdot T_2\cdots T_k\cdot [F]
$$
lies in $Z_0^{\mathrm{aff}}(\mathbb{R}^n)$. Since we work with tropical fan cycles, this is the origin equipped with an integer weight. We denote this integer by
$$
\deg(T_1\cdots T_k\cdot [F]).
$$

Finally, we conclude this section by recalling the projection formula, which will be used frequently throughout the paper.

\begin{definition}[Pushforward of cycles]\label{def:pushforward}
Let $\Phi\colon \mathbb R^m\to \mathbb R^n$ be a $\mathbb Z$-linear map, and let $[F]\in Z_k^{\mathrm{aff}}(\mathbb R^m)$. Choose a representative $(F,\omega_F)$ of $[F]$. After replacing $F$ by a suitable refinement, we may assume that the collection of images of those cones on which $\Phi$ is injective forms a rational fan in $\mathbb R^n$. The pushforward $\Phi_*[F]\in Z_k^{\mathrm{aff}}(\mathbb R^n)$ is represented by the weighted fan whose $k$-dimensional cones are the cones
$$
\sigma'=\Phi(\sigma),
\qquad
\sigma\in F^{(k)},\quad \dim \Phi(\sigma)=k.
$$
The weight of such a cone $\sigma'$ is defined by
$$
\omega_{\Phi_*F}(\sigma')
=
\sum_{\substack{\sigma\in F^{(k)}\\ \Phi(\sigma)=\sigma'}}
\omega_F(\sigma)\,
\bigl[\mathbb Z_{\sigma'}:\Phi(\mathbb Z_\sigma)\bigr],
$$
where $\mathbb Z_\sigma:=\mathbb Z^m\cap \operatorname{span}(\sigma)$ and $\mathbb Z_{\sigma'}:=\mathbb Z^n\cap \operatorname{span}(\sigma')$. Facets of $F$ on which $\Phi$ drops dimension do not contribute. The resulting affine cycle is independent of the chosen representative of $[F]$; cf. \cite[Construction~4.2, Definition~4.5, and
Proposition~4.6]{AllermannRau10}.
\end{definition}

Let $T:=\max(0,\ell_1,\dots,\ell_k)$ be a tropical regular function on $\mathbb R^n$. Its pullback along $\Phi$ is
the tropical regular function
$$
\Phi^*T\colon=T\circ \Phi=\max(0,\ell_1\circ\Phi,\dots,\ell_k\circ\Phi)
$$
on $\mathbb R^m$. Pushforward of cycles and pullback of
tropical regular functions are related by the following projection formula.

\begin{theorem}[{\cite[Projection formula]{AllermannRau10}}]
   Let $[C]\in Z^{\mathrm{aff}}_k(\mathbb{R}^n)$ be an affine $k$-cycle and let $\Phi:\mathbb{R}^n\to\mathbb{R}^m$ be an $\mathbb{Z}$-linear map such that $\Phi_*[C]\in Z^{\mathrm{aff}}_k(\mathbb{R}^m)$. Let $T$ be a tropical regular function over $\mathbb{R}^m$. Then the following holds:
   $$T\cdot\Phi_*[C]=\Phi_*(\Phi^*T\cdot [C]).$$
\end{theorem}

\begin{corollary}
    Let $[C]$ be an affine 1-cycle and let $\Phi:\mathbb{R}^n\to\mathbb{R}^m$ be an $\mathbb{Z}$-linear map such that $\Phi_*[C]$ is also a 1-cycle. Let $T$ be a tropical regular function over $\mathbb{R}^m$. Then the following holds:
    $$\deg(\Phi_*(\Phi^*T\cdot [C]))=\deg(\Phi^*T\cdot [C]).$$
\end{corollary}

We also recall that Bergman fans are tropical fans associated to loopless matroids. We use standard terminology from matroid theory, including flats, proper flats, bases, rank, and loopless matroids; see \cite{Oxley,ArdilaKlivans06,MaclaganSturmfels15}. For a realizable loopless matroid $M$, its Bergman fan $\mathcal{B}(M)$ may be described either by chains of proper flats or as the tropicalization of the corresponding constant-coefficient linear space; cf. \cite{ArdilaKlivans06,MaclaganSturmfels15}. Each maximal cone is given weight $1$.

\begin{definition}[Bergman fans and cycles]
    Let $M=(E,\mathcal{B})$ be a loopless matroid of rank r$(M)$ with ground set $E=\{0,\dots,n\}$. The Bergman fan $\mathcal{B}(M)$ is the fan
    $$\mathcal{B}(M):=\{w\in\mathbb{R}^n\mid\ M_w\ \rm{loopless}\}$$
    where $M_w$ is the associated matroid whose bases are the $w$-minimum bases of $M$.

    A tropical $k$-cycle $[F]$ is a Bergman cycle if $F=\mathcal{B}(M)$ for some loopless matroid $M$.
\end{definition}

\section{Tropical Regularity}\label{S3}

In this section, we introduce the notion of tropical regularity for tropical cycles and provide a conjectural classification of a large class of tropical regular cycles in terms of their combinatorial geometry.

\begin{definition}[Tropical regular cycles and regular tuples]
    A tropical $k$-cycle $[F]$ is said to be \textit{regular} if there exist $k$ tropical regular functions $T_1,\dots,T_k$, not necessarily distinct, such that
    $$
    \deg(T_1\cdots T_k\cdot [F])=1.
    $$
   We also call the tuple $(T_1,\dots,T_k,[F])$ regular if $\deg(T_1\dots T_k\cdot[F])=1$.
\end{definition}
We now introduce the key class of tropical (regular) cycles studied in this paper.

\begin{definition}\label{Hodgeindex}
    A tropical $k$-cycle $[F]$ satisfies the \emph{Hodge index inequality} if for every integer $p$ with $1 \le p \le k$ and all tropical regular functions $M_1,\dots,M_p,T_1,\dots,T_{k-p}$, we have
   $$
        \deg\bigl(M_1\cdots M_p \cdot T_1\cdots T_{k-p} \cdot [F]\bigr)^p
        \ge
        \prod_{i=1}^{p}
        \deg\bigl(\underbrace{M_i\cdots M_i}_{p\text{ times}} \cdot T_1\cdots T_{k-p} \cdot [F]\bigr).
  $$
\end{definition}

\begin{remark}
    This definition is parallel with the classical inequalities of Hodge
    type for nef divisor classes; see
    \cite[Section~1.6]{Laz04}.
    On a smooth projective surface, the Hodge index theorem gives
    $(A\cdot B)^2 \ge (A^2)(B^2)$ for any two nef divisor classes
    $A$ and $B$. Its higher-dimensional analogue gives the
    same mixed intersection inequality as in our definition. In particular, when $k=2$, the case $p=1$ is an equality. Thus a tropical
    $2$-cycle $[F]$ satisfies the Hodge index inequality if and only
    if
    $$
        \deg\bigl(T_1\cdot T_2\cdot[F]\bigr)^2
        \ge
        \deg\bigl(T_1\cdot T_1\cdot[F]\bigr)
        \deg\bigl(T_2\cdot T_2\cdot[F]\bigr)
    $$
    for all tropical regular functions $T_1$ and $T_2$.
\end{remark}

\begin{remark}
    The tropicalization $\operatorname{Trop}(U)$ of any irreducible
    subvariety $U\subseteq(\mathbb{C}^*)^n$ over the trivial valuation
    satisfies the Hodge index inequality.
    Indeed, its tropical intersection numbers agree with intersection
    numbers of nef divisors on a suitable projective compactification
    of $U$, where the classical Hodge inequalities apply
    \cite{Katz12,Laz04}.
    This class also contains Bergman fans of arbitrary loopless matroids
    \cite{MaclaganSturmfels15,AHK18}.
    These examples show that the definition describes a broad and
    natural class of tropical cycles, with sources in algebraic
    geometry, convex geometry, and matroid theory.
\end{remark}

\subsection{Tropical regular cycles satisfying the Hodge index inequality}

We now formulate a conjecture on the combinatorial structure of
tropical regular cycles satisfying the Hodge index inequality.
In Section~\ref{S5}, we prove this conjecture for tropical
regular $2$-fans.
We first introduce the notion of a tropical regular fibration:

\begin{definition}[Regular fibration]\label{fibration}
A regular tuple $(T_1,\dots,T_k,[F])$ forms a regular fibration if there exists a $\mathbb Z$-linear projection $p: \mathbb R^n\to \mathbb R^m$ such that the following conditions hold:
\begin{enumerate}
    \item The set-theoretic image $p(|F|)$ is the support of a regular tropical $r$-fan $G$ in $\R^m$, where $0<r<k$.
    
    \item There exists a nonempty subset
    $$
    S=\{s_1,\dots,s_r\}\subseteq \{1,\dots,k\}
    $$
    such that
    $$
    T_{s_i}=p^*(M_{s_i})
    \qquad \text{for all } s_i\in S,
    $$
    where $M_{s_i}$ are tropical regular functions on $\mathbb R^m$ satisfying
    $$
    \deg(M_{s_1}\cdots M_{s_r}\cdot [G])=1.
    $$
\end{enumerate}
We also call that a regular tuple forms a regular fibration along a $\Z$-linear projection $p$ if $p$ is clear.
\end{definition}

In our previous work, we gave a combinatorial classification of
tropical regular $1$-cycles \cite{Li25}.
Building on this result, Alexander Esterov raised the following
conjectural classification of higher-dimensional tropical regular
cycles satisfying the Hodge index inequality:

\begin{conjecture*}[A.Esterov]\label{trop}
     There exists a finite class $\mathcal{MD}_n$ of {\it minimal model} $k$-cycles, such that every regular $k$-cycle $[F]$ in $\R^n$ satisfying the Hodge index inequality and $\mathrm{span}_{\R}|F|=\R^n$, is one of the following: for any regular tuple $(T_1,\dots,T_k,[F])$,

\begin{itemize}
    \item[\textbf{I}:] either the regular tuple $(T_1,\dots,T_k,[F])$ forms a regular fibration,
    
    \item[\textbf{II}:] or there exists a saturated $\Z$-linear projection $\pi\colon\R^n\twoheadrightarrow\R^m$ such that $\pi_*[F]\in\mathcal{MD}_n$. Moreover, there exists a regular tuple $(M_1,\dots,M_k,\pi_*[F])$ such that
    $$T_i=\pi^*M_i$$
    for all $1\le i\le k$.\\
\end{itemize}
\end{conjecture*}

\begin{remark}
The notation \(\mathcal{MD}_n\) does not mean that all minimal-model
\(k\)-cycles in this finite class lie in the same ambient space \(\mathbb R^n\). In fact, \(\mathcal{MD}_n\) is a finite collection of \(k\)-cycles, possibly in different ambient spaces.
\end{remark}

\begin{remark}

A reader might hope for a stronger version of the conjecture: when $k>1$ and under the standing hypotheses, for every regular tropical cycle $[F]$ there is a single $\mathbb{Z}$-linear projection $\pi$ that works for all regular tuples, in the sense that either $(T_{1},\dots,T_{k},[F])$ forms a regular fibration along $\pi$, or $\pi_{*}[F]\in\mathcal{MD}_{n}$ together with $T_{i}=\pi^{*}M_{i}$ for all $1\le i\le k$. This stronger statement is false: a single projection cannot dominate all of the $T_{i}$ at once. We illustrate this with the following example.
\end{remark}

We close this section by introducing a useful theorem will be applied later:

\begin{theorem}[{Tropical Bernstein theorem \cite[Theorem 4.6.8]{MaclaganSturmfels15}}]\label{Bernstein} 
    Consider lattice polytopes $P_1,\dots,P_n$ in $\mathbb{R}^n$ with hypersurfaces $\Sigma_{P_1},\dots,\Sigma_{P_n}$. Then the stable intersection 
    $$\Sigma_{P_1}\cap_{st}\dots\cap_{st}\Sigma_{P_n}$$
    is the origin $0$ with weight given by the mixed volume $\mathrm{MV}(P_1,\dots,P_n)$.
\end{theorem}

\begin{example}
Let $P=\mathrm{conv}\{0,(1,0,0),(0,1,0),(1,1,0),(1,0,5),(0,1,6),(0,0,7),(1,1,8)\}\subset\mathbb{R}^{3}$ and let $[F]:=[\Sigma_{P}]$ be its associated hypersurface in $(\mathbb{R}^{3})^{\vee}$; since $\dim P=3$, the cycle $[F]$ spans $(\mathbb{R}^{3})^{\vee}$. Set $$T_{1}=\max(0,\,x,\,x+z),\qquad T_{2}=\max(0,\,y,\,y+z),$$ with Newton polytopes the unimodular triangles $Q_{1}=\mathrm{conv}\{0,(1,0,0),(1,0,1)\}$ and $Q_{2}=\mathrm{conv}\{0,(0,1,0),(0,1,1)\}$. The Esterov--Gusev theorem (Theorem~\ref{E}) therefore gives $\mathrm{MV}(Q_{1},Q_{1},P)=1$, and Theorem \ref{Bernstein}  yields $\deg(T_{1}\cdot T_{1}\cdot[F])=1$. Similarly, we obtain $\deg(T_{2}\cdot T_{2}\cdot[F])=1$ as well.

Suppose there is a $\mathbb{Z}$-linear projection $\pi\colon(\mathbb{R}^{3})^{\vee}\to(\mathbb{R}^{2})^{\vee}$ such that $T_{1}=\pi^{*}M_{1}$ and $T_{2}=\pi^{*}M_{2}$. Since $\pi^{*}M_{i}=M_{i}\circ\pi$ is constant along $\ker\pi$, and a tropical regular function $\max(0,\ell_{1},\dots,\ell_{q})$ is constant along a vector precisely when that vector lies in $\bigcap_{j}\ker \ell_{j}=L(T_{i})$, we would have $\ker\pi\subseteq L(T_{1})\cap L(T_{2})$. However, we have$$L(T_{1})=\{v:v_{1}=v_{3}=0\}=\mathbb{R}e_{2},\qquad L(T_{2})=\{v:v_{2}=v_{3}=0\}=\mathbb{R}e_{1},$$ so $L(T_{1})\cap L(T_{2})=\{0\}$. Hence $\ker\pi=\{0\}$, forcing $\pi$ to be injective, contradiction. Therefore no single projection dominates both $T_{1}$ and $T_{2}$. \\
\end{example}

\section{Tropical regular $2$-cycles} \label{S4} 

In this section, we aim to show that conjecture \hyperref[trop]{A.Esterov} holds for regular tropical $2$-cycles satisfying the Hodge index inequality:

\begin{theorem}\label{Thm_1}
     There exists a finite class $\mathcal{MD}_n$ of {\it minimal model} $2$-cycles, such that every regular $2$-cycle $[F]$ in $\R^n$ satisfying the Hodge index inequality and $\mathrm{span}_{\R}|F|=\R^n$, is one of the following: for any regular tuple $(T_1,T_2,[F])$,

\begin{itemize}
    \item[\textbf{I}:] either the regular tuple $(T_1,T_2,[F])$ forms a regular fibration,
    
    \item[\textbf{II}:] or there exists a saturated $\Z$-linear projection $\pi\colon\R^n\twoheadrightarrow\R^m$ such that $\pi_*[F]\in\mathcal{MD}_n$. Moreover, there exists a regular tuple $(M_1,M_2,\pi_*[F])$ such that
    $$T_i=\pi^*M_i$$
    for all $i=1,2$.
\end{itemize}
\end{theorem}

Recall that $\mathrm{TR}(\R^n)$ denotes the family of tropical regular
functions on $\R^n$ (cf. Remark~\ref{TR(Rn)}).
We now extend this notion by allowing rational coefficients.
A function $T:\R^n\to\R$ is called a
\textit{quasi-tropical regular function} if $ T=\max\{0,q_1,\dots,q_r\}$ for some $r\in\Z_{>0}$ and $q_1,\dots,q_r\in(\Q^n)^\vee$. We denote the family of these functions by $\mathrm{QTR}(\R^n)$. For any $T=\max\{0,q_1,\dots,q_r\}\in\mathrm{QTR}(\R^n)$,
clearing denominators gives a positive integer $\delta$ such that $ \delta T=\max\{0,\delta q_1,\dots,\delta q_r\}
    \in\mathrm{TR}(\R^n)$. Clearly, $\mathrm{TR}(\R^n)\subset\mathrm{QTR}(\R^n)$. 
Moreover, $T_1+T_2\in\mathrm{QTR}(\R^n)$ for any $T_1,T_2\in\mathrm{QTR}(\R^n)$ by taking the Minkowski sum of their Newton polytopes.

\begin{definition}
    A tropical $2$-cycle $[F]$ in $\R^n$ determines a bilinear form $$Q_{[F]}:\mathrm{QTR}(\R^n)\times\mathrm{QTR}(\R^n)\to\Q_{\ge 0};\qquad (T_1,T_2)\mapsto \frac{1}{\delta\zeta}\deg(\zeta T_1\cdot \delta T_2\cdot [F])$$
    for any $\delta,\zeta\in\Z_{>0}$ makes $\zeta T_1,\delta T_2\in \mathrm{TR}(\R^n)$.
\end{definition}

This definition is well-defined: assume $T_2\in\mathrm{TR}(\R^n)$ while $T_1\in\mathrm{QTR}(\R^n)$. If $T_1\in\mathrm{TR}(\R^n)$ then it follows from (4) of proposition \ref{propTCI} directly. If $T_1\notin \mathrm{TR}(\R^n)$ then there is a minimal positive integer $\zeta$ makes $\zeta T_1\in\mathrm{TR}(\R^n)$ while any other $\zeta'\in\{\zeta'\in\Z_{>0}:\zeta' T_1\in\mathrm{ TR}(\R^n)\}$ can be divided by $\zeta$. Therefore
$$Q_{[F]}(T_1,T_2)=\frac{1}{\zeta}\deg(\zeta T_1\cdot T_2\cdot [F])=\frac{\alpha}{\alpha\zeta}\deg(\zeta T_1\cdot T_2\cdot [F])=\frac{1}{\alpha\zeta}\deg(\alpha\zeta T_1\cdot T_2\cdot [F])$$
for any $\alpha\in\Z_{>0}$. This shows that the definition is well-defined. It is clear that $Q_{[F]}(\cdot,\cdot)$ is also commutative. Moreover, if $[F]$ satisfies the Hodge index inequality (cf. Definition \ref{Hodgeindex}), then $$Q_{[F]}(A,B)^2\ge Q_{[F]}(A,A)Q_{[F]}(B,B)$$for any $A,B\in\mathrm{QTR}(\R^n)$.

\begin{proposition}
    $Q_{[F]}((\epsilon T_1+T_2),T_3)=\epsilon Q_{[F]}(T_1,T_3)+Q_{[F]}(T_2,T_3)$ for any $T_1,T_2,T_3\in\mathrm{QTR}(\R^n)$ and $\epsilon\in\Q_{\ge 0}$.
\end{proposition}

\begin{proof}
    Clearly $\epsilon Q_{[F]}(T_1,T_3)=Q_{[F]}(\epsilon T_1,T_3)$, it is enough to show $Q_{[F]}((T_1+T_2),T_3)=Q_{[F]}(T_1,T_3)+Q_{[F]}(T_2,T_3)$. Assume $T_3\in\mathrm{TR}(\R^n)$, $\zeta T_1,\delta T_2\in\mathrm{ TR}(\R^n)$. Then $Q_{[F]}((T_1+T_2),T_3)$ is 
    $$\frac{1}{\delta\zeta}\deg(\delta\zeta(T_1+T_2)\cdot T_3\cdot [F])=\frac{1}{\delta\zeta}\left( \deg(\delta\zeta T_1\cdot T_3\cdot [F])+\deg(\delta\zeta T_2\cdot T_3\cdot [F])\right)=Q_{[F]}(T_1,T_3)+Q_{[F]}(T_2,T_3)$$as required.\\
\end{proof}

\subsection{Toolkit} We will state and prove the lemmas which will be used in the proof of Theorem \ref{Thm_1}.

\begin{lemma}\label{4L1}
    Let $[F]$ be a nonzero tropical $2$-cycle in $\R^n$. There is a representative $F$ of $[F]$ and a tropical regular function $K\in\mathrm{TR}(\R^n)$ such that $K\cdot[F]$ is a $1$-cycle and
    $$\{\rho\in K\cdot [F]\}=\{\rho\in F(1)\}.$$
    In other words, the rays in the $1$-cycle $K\cdot [F]$ are exactly the rays in $F$. Moreover $\deg(K\cdot K\cdot [F])>0$.
\end{lemma}

\begin{proof}
 We begin with an arbitrary representative $F'$ of $[F]$ and describe each of its $2$-cones by rational homogeneous linear equalities and inequalities. Clearing denominators gives nonzero integral linear functionals $\ell_1,\dots,\ell_s$. The hyperplanes $\ker(\ell_1),\dots,\ker(\ell_s)$ form a central hyperplane arrangement in $\R^n$.
We denote by $F$ the subdivision of $F'$ induced by this arrangement. Define 
$$K(x)=\sum_{i=1}^s|\ell_i(x)|.$$
This is a globally integral piecewise linear function in $\mathrm{TR}(\R^n)$. Indeed, for any $\epsilon=(\epsilon_1,\dots,\epsilon_s)\in\{-1,0,1\}^s$, we define a linear functional $h_\epsilon:=\sum_{i=1}^s\epsilon_i\ell_i$. Therefore $K(x)=\max_{\epsilon\in\{-1,0,1\}^s}(0,h_\epsilon)$ which is exactly a tropical regular function. Note that $K(x)$ is linear on every chamber of the hyperplane arrangement, so it is linear on every $2$-cone of $F$. 

Let $\tau$ be a ray of $F$ with primitive generator $u_\tau$, the restriction of $K$ on $\tau$ is the linear form
$$K_\tau=:\lambda_\tau=\sum_{i=1}^s\mathrm{sign}(\ell_i(u_\tau))\ell_i;\qquad\mathrm{sign}(0)=0.$$
Therefore, we have$$K-\lambda_\tau=\sum_{i=1}^s|\ell_i|-\mathrm{sign(\ell_i(u_\tau))}\ell_i.$$
Since each $|\ell_i|-\mathrm{sign}(\ell_i(u_\tau))\ell_i$ is non-negative, the kernel of $K-\lambda_\tau$ is the set of vectors satisfying the following:
$$\begin{cases}
    \ell_i(v)\ge 0&\text{if }\ell_i(u_\tau)>0\\
    \ell_i(v)\le 0&\text{if }\ell_i(u_\tau)<0\\
    \ell_i(v)=0 &\text{if }\ell_i(u_\tau)=0.
\end{cases}
$$ Hence we conclude that
\begin{equation}\{v\in F\mid K(v)=\lambda_\tau(v)\}=\tau.\end{equation}
Note that for any $2$-cone $\sigma\subset F$ contains $\tau$ as a face, there always exists a primitive generator $v_{\sigma/\tau}$ of $\Z_\sigma/\Z_\tau$ lives in $\sigma\setminus\tau$ which can be obtained by starting from arbitrary primitive generator $w_{\sigma/\tau}$ and add a sufficiently large integral multiple of $u_\tau$. Now the formula of tropical intersection product (definition \ref{TCI}) gives the weight of $\tau$ in $K\cdot [F]$:
    \begin{align*}
\omega_{K\cdot[F]}(\tau)&=\sum_{\sigma\in F^{(2)},\tau<\sigma}K(\omega_F(\sigma)v_{\sigma/\tau})-K_\tau\left(\sum_{\sigma\in F^{(2)},\tau<\sigma}\omega_F(\sigma)v_{\sigma/\tau}\right)\\
&=\sum_{\sigma\in F^{(2)},\tau<\sigma}\omega_F(\sigma)\left(K(v_{\sigma/\tau})-\lambda_\tau(v_{\sigma/\tau)} \right).
    \end{align*}
The equation (4.1) yields that $K(v)-\lambda_\tau(v)=(K-\lambda_\tau)(v)>0$ for all $v\notin\tau$. Therefore $\omega_{K\cdot[F]}(\tau)>0$ for every ray $\tau$ of $F$ since $\omega_F(\sigma)>0$ and $K(v_{\sigma/\tau})-\lambda_\tau(v_{\sigma/\tau})>0$. The latter assertion holds immediately.\\
\end{proof}

\begin{lemma}\label{4L2}
    Two tropical $1$-cycles $C_1,C_2$ in $\R^n$ are equal if $$\deg(D\cdot C_1)=\deg(D\cdot C_2)$$
    for all $D\in\mathrm{TR}(\R^n)$.
\end{lemma}

\begin{proof}
    The case $n=1$ follows by taking $D(x)=\max\{0,x\}$.

    For $n=2$, fix a ray $\rho$ occurring in either $C_1$ or $C_2$ with primitive
generator $u_\rho$, and let $\omega_i(\rho)$ be its weight in $C_i$ while $\omega_i(\rho)=0$ if $\rho\notin C_i$.
Choose integral linear forms $\ell,m$ such that
$\ker\ell=\R u_\rho$ and $m(u_\rho)>0$, and set
$$
    U=\max\{0,N\ell,-N\ell\},
    \qquad
    V=\max\{0,N\ell,-N\ell,m\}.
$$
For sufficiently large positive integer $N$ we have $U(u_\tau)=V(u_\tau)$ for any other primitive generators $u_\tau$ of the ray in either $C_1$ or $C_2$ while both functions vanish at $-u_\rho$.
Since $V(u_\rho)-U(u_\rho)=m(u_\rho)$, so we have
$$
    \deg(V\cdot C_i)-\deg(U\cdot C_i)
    =m(u_\rho)\omega_i(\rho).
$$
Applying the hypothesis separately to $U$ and $V$ yields 
$$\deg(V\cdot C_1)-\deg(U\cdot C_1)=m(u_\rho)\omega_1(\rho)=m(u_\rho)\omega_2(\rho)=\deg(V\cdot C_2)-\deg(U\cdot C_2)$$
while the positivity of $m(u_\rho)$ forces $\omega_1(\rho)=\omega_2(\rho)$. As $\rho$ was arbitrary, we obtain $C_1=C_2$.

    Let $k>2$, and assume that the statement holds for all $2\le n<k$.
    Suppose that $C_1,C_2$ are tropical $1$-cycles in $\R^k$ such that $\deg(D\cdot C_1)=\deg(D\cdot C_2)$ for all $D\in\mathrm{TR}(\R^k)$. Now we define $S=C_1\cup C_2$, and put $\langle S\rangle:=\{\operatorname{span}(\rho_1,\rho_2) \mid \rho_1,\rho_2\in S\}.$
    We can choose a rational line $L$ not contained in any
    $2$-plane of $\langle S\rangle$.
    Then choose a $\Z$-linear projection
    $\pi:\R^k\to\R^{k-1}$ with kernel $L$.

    The restriction of $\pi$ on $S$ is injective: if $x\in\rho_1$, $y\in\rho_2$, and $\pi(x)=\pi(y)$, then
      $ x-y\in L\cap\operatorname{span}(\rho_1,\rho_2)=\{0\}$. Therefore $\pi_*C_1\neq\pi_*C_2$. Indeed, if $C_1$ and $C_2$ have the same support but different weights,
the two weights on each ray are multiplied by the same positive
lattice index under $\pi$, so the pushforward cycles remain distinct.
If their supports differ, injectivity ensures that their image
supports also differ.
    For every $M\in\mathrm{TR}(\R^{k-1})$, the pullback $\pi^*M\in\mathrm{TR}(\R^k)$.
    The hypothesis therefore gives
   $$
        \deg(\pi^*M\cdot C_1)=\deg(\pi^*M\cdot C_2).
    $$
    By the projection formula,
    $$
        \deg(M\cdot\pi_*C_1)=\deg(M\cdot\pi_*C_2).
   $$
    The induction hypothesis now implies
    $\pi_*C_1=\pi_*C_2$, a contradiction.\\
\end{proof}

\begin{lemma}\label{4L3}
    Let $[F]$ be a tropical $2$-cycle in $\R^n$ such that $\mathrm{span}_{\R}|F|=\R^n$ and satisfies the Hodge index inequality and let $T\in\mathrm{TR}(\R^n)$ such that $Q_{[F]}(T,T)>0$. Then the tropical $1$-cycle $T\cdot [F]$ spans $\R^n$.
\end{lemma}

\begin{proof}
    Consider $[C]:=T\cdot [F]$ and $\alpha:=Q_{[F]}(T,T)$. If $\mathrm{span}_\R|C|\neq\R^n$, then there exists a $0\neq\ell\in(\Z^n)^\vee$ such that $\ell(x)=0$ for all $x\in\mathrm{span}_\R|C|$. Therefore $\deg(\max(0,\ell)\cdot T\cdot [F])=\deg(\max(0,\ell)\cdot[C])=0$. Moreover, $\deg(\max(0,\ell)\cdot\max(0,\ell)\cdot[F])=0$ since $\max(0,\ell)\cdot[F]$ is a $1$-cycle supported in $\ker\ell$. Consequently we have $Q_{[F]}(\max(0,\ell),\max(0,\ell))=0$ and $Q_{[F]}(\max(0,\ell),T)=0$. 

    Now for any $M\in\mathrm{TR}(\R^n)$ and any $k\in\Z_{>0}$, the Hodge index inequality gives
        \begin{align*}
Q_{[F]}(T,M+k\max(0,\ell))^2&\ge Q_{[F]}(T,T)Q_{[F]}(M+k\max(0,\ell),M+k\max(0,\ell))\\
            &\ge Q_{[F]}(T,T)\left(Q_{[F]}(M,M)+k^2 Q_{[F]}(\max(0,\ell),\max(0,\ell))+2kQ_{[F]}(M,\max(0,\ell)) \right)\\
            &\ge Q_{[F]}(T,T)\left(Q_{[F]}(M,M)+2kQ_{[F]}(M,\max(0,\ell))\right)
        \end{align*}
while the last inequality comes from $Q_{[F]}(\max(0,\ell),\max(0,\ell))=0$. The left hand side has
\begin{align*}
    Q_{[F]}(T,M+k\max(0,\ell))^2&=\left(Q_{[F]}(T,M)+kQ_{[F]}(T,\max(0,\ell))\right)^2\\
    &=Q_{[F]}(T,M)^2
\end{align*}
while the last equality comes from $Q_{[F]}(T,\max(0,\ell))=0$. Therefore, we conclude that
$$Q_{[F]}(T,M)^2\ge Q_{[F]}(T,T)\left(Q_{[F]}(M,M)+2kQ_{[F]}(M,\max(0,\ell))\right).$$
The condition $Q_{[F]}(T,T)>0$ forces $Q_{[F]}(M,\max(0,\ell))=0$ since $k$ can be arbitrarily large. In other words, $\max(0,\ell)\cdot[F]$ is a tropical $1$-cycle $[D]$ such that $\deg(M\cdot [D])=0$ for every $M\in\mathrm{TR}(\R^n)$. Finally, $|F|$ spanning $\R^n$ forces $[D]\neq 0$, hence there must exists an $M\in\mathrm{TR}(\R^n)$ such that $\deg(M\cdot\max(0,\ell)\cdot [F])>0$, contradiction. Therefore, we conclude that $T\cdot [F]=[C]$ spans $\R^n$ as required.\\
\end{proof}

\begin{theorem}[\cite{Li25} Theorem 4.13]
\label{4L4}
    Assume that $T_1,T_2$ are two tropical regular functions over $\mathbb{R}^n$ such that the Minkowski summation has dim$(\mathrm{Newt}(T_1)+\mathrm{Newt}(T_2))=n$. Then there are finitely many tropical $2$-cycles $[F]$ not necessarily spanning $\R^n$ such that $\deg(T_i\cdot T_j\cdot [F])\le 1$ for all $i,j\in\{1,2\}$.\\
\end{theorem}

\subsection{Proof}

In this subsection, we will prove Theorem~\ref{Thm_1}.
Recall that $[F]$ is a tropical regular $2$-cycle in $\R^n$
whose support spans $\R^n$ satisfying the Hodge index
inequality.
Thus, for all $T_1,T_2\in\mathrm{TR}(\R^n)$, we have
\[
    Q_{[F]}(T_1,T_2)^2
    \ge
    Q_{[F]}(T_1,T_1)\,Q_{[F]}(T_2,T_2).
\]

We fix an arbitrary regular tuple $(T_1,T_2,[F])$, so that
$Q_{[F]}(T_1,T_2)=1$.
We show that if $Q_{[F]}(T_i,T_i)=0$ for some $i\in\{1,2\}$,
then we obtain the regular fibration case of
Theorem~\ref{Thm_1}.
If instead
\[
    Q_{[F]}(T_1,T_1)=Q_{[F]}(T_2,T_2)=1,
\]
we obtain the second case of the theorem. \\

\begin{proposition}\label{4p1}
    Let $[F]$ be a tropical regular $2$-cycle in $\R^n$ such that $\mathrm{span}_\R|F|=\R^n$ and satisfying the Hodge index inequality. If there exist $T_1,T_2\in\mathrm{TR}(\R^n)$ such that $\deg(T_1\cdot T_2\cdot [F])=1$ and $\deg(T_1\cdot T_1\cdot [F])=0$, then the tuple $(T_1,T_2,[F])$ forms a regular fibration.
\end{proposition}

\begin{proof}
    Consider the $\Z$-linear projection $\pi:\R^n\to\R^n/L(T_1)\cong\R^m$ where $L(T_1)$ is the lineality space of $T_1$ (cf. Definition \ref{TR}), then $T_1=\pi^*M$ for some $M\in\mathrm{TR}(\R^m)$ such that $\mathrm{dim}(\mathrm{Newt}(M))=m$. The projection formula yields
    $$\deg(T_1\cdot T_1\cdot [F])=0=\deg(M\cdot M\cdot \pi_*[F]).$$
    First, $|\pi(F)|\neq\{0\}$, since $|F|$ spans $\R^n$ and is
therefore not contained in $L(T_1)$. 
If $|\pi(F)|$ is a $2$-fan, then $\pi_*[F]$ would be a
nonzero tropical $2$-cycle and this forces $\deg(M\cdot M\cdot\pi_*[F])>0$ since $\mathrm{Newt}(M)$ is full-dimensional. Therefore, $\pi(F)$ is a $1$-fan.

Applying Lemma \ref{4L1} to $[F]$ we obtain a function $K$ and a refinement $F$. On the one hand, $|\pi_*(K\cdot[F])|\subseteq|\pi(F)|$ because $|K\cdot [F]|\subset|F|$. On the other hand, any ray $\rho\subset\pi(F)$ corresponds to at least one ray $r\subset F$ such that $\pi(r)=\rho$ while $r$ is also a ray of $K\cdot [F]$ by Lemma \ref{4L1}. Therefore $r\subset|\pi_*(K\cdot [F])|$, i.e., $|\pi(F)|\subseteq|\pi_*(K\cdot [F])|$. Combine them together we have
\begin{equation}|\pi_*(K\cdot [F])|=|\pi(F)|.\end{equation}
We aim to show that $|\pi(F)|=|\pi_*(T_2\cdot [F])|$ (\textbf{Target} $\mathcal{A}$), then the projection formula implies that $\pi_*(T_2\cdot [F])$ is regular:
$$\deg(M\cdot\pi_*(T_2\cdot [F]))=\deg(\pi^*M\cdot T_2\cdot [F])=\deg(T_1\cdot T_2\cdot [F])=1.$$
Thus $(T_1,T_2,[F])$ forms a regular fibration along $\pi$ as required.

Consider the $1$-cycle $\pi_*((K+T_2)\cdot [F])$, equation (4.2) deduces that \begin{equation}|\pi_*((K+T_2)\cdot[F])|=|\pi(F)|.\end{equation} Now for any $D\in\mathrm{TR}(\R^m)$, $\pi_*[F]=0$ gives
\begin{equation}Q_{[F]}(T_1,\pi^*D)=0=Q_{[F]}(\pi^*D,\pi^*D)\end{equation}
via the projection formula. Write $a=Q_{[F]}(T_1,K+T_2)$, $b=Q_{[F]}(T_2,\pi^*D)$ and $c=Q_{[F]}(K+T_2,\pi^*D)$. We will show that $c=ab$ (\textbf{Target} $\mathcal{B}$), so that 
$$\deg(D\cdot\pi_*((K+T_2)\cdot[F]))=\deg(D\cdot a\left(\pi_*(T_2\cdot[F])\right))$$ holds for any $D\in\mathrm{TR}(\R^m)$.
Lemma \ref{4L2} yields $a(\pi_*(T_2\cdot[F]))=\pi_*((K+T_2)\cdot[F])$ whence the \textbf{Target} $\mathcal{A}$ holds: $|\pi(F)|\underset{(4.3)}{=}|\pi_*((K+T_2)\cdot[F])|=|\pi_*(T_2\cdot [F])|$. So we can solve \textbf{Target} $\mathcal{A}$ by proving \textbf{Target} $\mathcal{B}$.

To see $c=ab$, we introduce two positive rational variables $r,t\in\Q_{>0}$. Consider two quasi-tropical regular functions $M_1:=T_1+tT_2$ and $M_2:=\pi^*D+rt(K+T_2)$:
\begin{align*}
    Q_{[F]}(M_1,M_2)&=Q_{[F]}(T_1,\pi^*D)+Q_{[F]}(T_1,rt(K+T_2))+Q_{[F]}(tT_2,\pi^*D)+Q_{[F]}(tT_2,rt(K+T_2))\\
    &\underset{(4.4)}{=}0+rt\cdot a+t\cdot b+rt^2\cdot Q_{[F]}(T_2,(K+T_2));
\end{align*}
\begin{align*}
    Q_{[F]}(M_1,M_1)&=Q_{[F]}(T_1,T_1)+2tQ_{[F]}(T_1,T_2)+t^2Q_{[F]}(T_2,T_2)\\
    &=0+2t+t^2Q_{[F]}(T_2,T_2);
\end{align*}
\begin{align*}
    Q_{[F]}(M_2,M_2)&=Q_{[F]}(\pi^*D,\pi^*D)+2rtQ_{[F]}(\pi^*D,K+T_2)+r^2t^2Q_{[F]}(K+T_2,K+T_2)\\
    &\underset{(4.4)}{=}0+2rt\cdot c+r^2t^2Q_{[F]}(K+T_2,K+T_2).
\end{align*}
Recall that $[F]$ satisfies the Hodge index inequality, so that $Q_{[F]}(M_1,M_2)^2\ge Q_{[F]}(M_1,M_1)Q_{[F]}(M_2,M_2)$:
\begin{align*}
    \left(rta+tb+rt^2Q_{[F]}(T_2,K+T_2)\right)^2&\ge\left(2t+t^2Q_{[F]}(T_2,T_2)\right)\left(2rtc+r^2t^2Q_{[F]}(K+T_2,K+T_2)\right)
\end{align*}
After dividing both sides by $t^2$ (since $t>0$), we obtain
$$\left(ra+b+rtQ_{[F]}(T_2,K+T_2)\right)^2\ge\left(2+tQ_{[F]}(T_2,T_2)\right)\left(2rc+r^2tQ_{[F]}(K+T_2,K+T_2)\right).$$
Recall that $a,b,c$ are fixed non-negative integers, while the inequality above holds for arbitrary $r,t\in Q_{>0}$. Hence we fix $r$ and let $t\to 0$ through positive rational values. This gives
\begin{equation}(ra+b)^2\ge4rc.\end{equation}
Applying the same argument to $N_1:=T_1+t(K+T_2)$ and $N_2:=\pi^*D+rtT_2$ we obtain
\begin{equation}
    (c+r)^2\ge 4rab.
\end{equation}
Note that $a=Q_{[F]}(T_1,K+T_2)=Q_{[F]}(T_1,K)+Q_{[F]}(T_1,T_2)\ge 1$ but $b,c$ can be $0$.

If $b=0$, then equation (4.5) yields $\frac{ra^2}{4}\ge c$ which forces $c=0$ by letting $r\to 0$; if $b>0$, then we substitute $r$ by $b/a$ into (4.5) and obtain $ab\ge c$.

If $c=0$, then equation (4.6) yields $r\ge 4ab$ which forces $b=0$ by letting $r\to 0$; if $c>0$ then we substitute $r$ by $c$ in (4.6) and obtain $c\ge ab$.

Combine these together, we conclude that $c=ab$ so that the \textbf{Target} $\mathcal{B}$ follows as required and this finishes the proof.\\
\end{proof}

\begin{proposition}\label{4p2}
    Let $[F]$ be a tropical regular $2$-cycle in $\R^n$ such that $\mathrm{span}_\R|F|=\R^n$ and satisfying the Hodge index inequality. If there exist $T_1=\max(0,\ell_1,\dots,\ell_t),T_2=\max(0,\mu_1,\dots,\mu_s)\in\mathrm{TR}(\R^n)$ such that
    $$\deg(T_1\cdot T_1\cdot [F])=\deg(T_1\cdot T_2\cdot [F])=\deg(T_2\cdot T_2\cdot [F])=1,$$
    then there exists a saturated $\Z$-linear projection $\pi:\R^n\twoheadrightarrow\R^m$ ($2\le m\le n$) and $M_1,M_2\in\mathrm{TR}(\R^m)$ satisfying:
    \begin{enumerate}
        \item $T_1=\pi^*M_1$, $T_2=\pi^*M_2$;

        \item $\mathrm{Newt}(M_1),\mathrm{Newt}(M_2)\subseteq\mathrm{conv}\{0,\pm e_i,e_i-e_j\mid 1\le i\neq j\le m\}$;

        \item  $\mathrm{dim}(\mathrm{Newt}(M_1)+\mathrm{Newt}(M_2))=m$;

        \item $\deg(M_i\cdot M_j\cdot \pi_*[F])=1$ for $i,j\in\{1,2\}$. 
    \end{enumerate}
\end{proposition}

\begin{proof}
    Lemma \ref{4L3} deduces that $[C]=T_1\cdot [F]$ is a tropical $1$-cycle spanning $\R^n$. We denote the primitive generator of $[C]$ by $v_i$ with weights $\omega_i$. Note that $\deg(T_1\cdot [C])=1$ determines a unique index $a$ such that
    $$\omega_a=1,\qquad T_1(v_a)=1,\qquad T_1(v_i)=0\ \ (i\neq a).$$
    Moreover, any non-zero linear form $\ell_i$ of $T_1:=\max(0,\ell_1,\dots,\ell_t)$ satisfies the followings \cite{Li25}:
    \begin{enumerate}
        \item $\ell_i(v_a)=1$;

        \item there exists a unique primitive generator $v_i$ such that $\ell_i(v_i)=-1$ and $\ell_i(v_j)=0$ for all $j\neq a,i$.
    \end{enumerate}
    In other words, every nonzero linear form $\ell_i$ determines a gallery of $[C]$, cf. \cite{Li25}. 

    Applying the same argument to $\deg(T_2\cdot [C])=1$ yields a unique index $b$ for a primitive generator $v_b$, not necessarily distinct with $a$, and the same results on $\mu_i$'s. 

    We now define a directed graph inspired by the gallery construction for regular
$1$-cycles. Its vertex set is
\[
S := \bigl\{v_q \mid
\ell_i(v_q)\neq 0 \text{ for some } i
\text{ or }
\mu_j(v_q)\neq 0 \text{ for some } j
\bigr\},
\]
where the $v_q$ are the primitive ray generators, and the $\ell_i$ and
$\mu_j$ are the defining linear functionals of $T_1$ and $T_2$, respectively. By the discussion above, all nonzero $\ell_i$ take the value $1$ at
$v_a$, and all nonzero $\mu_j$ take the value $1$ at $v_b$, possibly with
$v_a=v_b$. For each nonzero $\ell_i$ (respectively, $\mu_j$), draw an edge
from $v_a$ (respectively, $v_b$) to the unique vertex where it takes the
value $-1$. Hence we obtain a directed graph $G$. If $G$ is connected, choose a spanning tree rooted at $v_a$. Otherwise, choose a spanning tree in each of its two connected components, rooted at $v_a$ and $v_b$, respectively. Since any edge $e$ of the tree is either $v_a\to w$ or $v_b\to w$ for some $w\in G$, we assign a linear form $\lambda_e$ of this edge $e$ such that $\lambda_e=\ell_w$ if $e:v_a\to w$ while $\lambda_e=\mu_w$ if $e:v_b\to w$.

For a component rooted $v_r$ ($r=a$ or $b$), set $h_r=0\in(\Z^n)^\vee$. For every other vertex
$v_j$ in that component, we associate a linear form $h_j$ to $v_j$:
\begin{equation}
h_j := \sum_{e\in P_{rj}}\varepsilon_{rj}(e)\lambda_e \in (\Z^n)^\vee,
\end{equation}
where $P_{rj}$ is the unique tree path from $v_r$ to $v_j$, and
$\varepsilon_{rj}(e)$ is $1$ if the path traverses $e$ along its chosen
orientation and $-1$ otherwise.

Note that $\lambda_e$ takes the values $1$ and $-1$ at its tail and head but vanishes on all other primitive generators. This provides an expression of $h_j$: for any primitive generator $v_i$
$$h_j(v_i)=\begin{cases}
    1,&\text{if }i=r,\\
    -1,&\textbf{if }i=j,\\
    0,&\text{otherwise}.
\end{cases}$$
Notice that we have exactly $m=|\mathrm{V}(G)|-c(G)$ where $\mathrm{V}(G)$ is the vertex set and $c(G)$ is the number of connected components, so we might re-index those non-root vertices by $v_1,\dots,v_m$. Consequently, these $h_j$'s are linearly independent in $(\R^n)^\vee$. In other words, $m\le n$. Now we define a $\Z$-linear map 
$$\pi:\R^n\to \R^m,\qquad x\mapsto (h_1(x),\dots,h_m(x)).$$
Since $\pi:v_i\mapsto -e_i$ where $\{e_1,\dots,e_m\}$ is the standard $\Z$-basis of $\Z^m$, we conclude that $\pi$ is a saturated $\Z$-projection. 

Since $|C|$ spans $\R^n$, linear forms are uniquely determined by their values on the primitive generators of $C$. As shown above, each linear form $\ell_i$ of $T_1$ determines a unique $v_i$ such that
$\ell_i(v_i)=-1$, $\ell_i(v_a)=1$, and $\ell_i(v_j)=0$
for all $j\notin\{i,a\}$.

If $G$ is disconnected, its connected components are themselves trees
by construction, so $\ell_i=h_i$.

If $G$ is connected and $v_a=v_b$, then $G$ itself is the spanning tree,
and again $\ell_i=h_i$.

If $G$ is connected and $v_a\neq v_b$, then $\ell_i=h_i$ when the root
is $v_a$, and $\ell_i=h_i-h_a$ when the root is $v_b$. In particular, if there exists an edge $v_a\to v_b$, then $\ell_b=-h_a$ when $v_b$ is the root.

Consequently, we have $\ell_i=\pi^*x_i$, $\ell_i=\pi^*(x_i-x_a)$ or $\ell_i=\pi^*(-x_i)$. The same argument holds for all $\mu_i$ of $T_2$ as well. Therefore, we finish the construction of $\pi$ and proof of (1) and (2).

Part (3) follows directly from equation (4.7): every $h_j$ is a linear combination of some $\ell$'s and $\mu$'s while the dimension of $\mathrm{span}(h_1,\dots,h_m)=m$. Hence $\mathrm{Newt}(M_1)+\mathrm{Newt}(M_2)$ has dimension $m$ as required.

The projection formula gives
$$\deg(M_i\cdot M_j\cdot \pi_*[F])=1\qquad(i,j\in\{1,2\})$$
which proves (4).\\

\end{proof}

\begin{proof}[Proof of Theorem \ref{Thm_1}]
Let $[F]$ be a tropical $2$-cycle in $\R^n$ satisfying the Hodge index
inequality and $\mathrm{span}_\R|F|=\R^n$, and let
$T_1,T_2\in\mathrm{TR}(\R^n)$ satisfy
$\deg(T_1\cdot T_2\cdot[F])=1$. The Hodge index inequality gives
\[
1=\deg(T_1\cdot T_2\cdot[F])^2
\ge
\deg(T_1\cdot T_1\cdot[F])\,
\deg(T_2\cdot T_2\cdot[F]).
\]
Since the self-intersection numbers are non-negative,
there are two possibilities:
\begin{enumerate}
    \item $\deg(T_i\cdot T_i\cdot[F])=0$ for some $i\in\{1,2\}$;
    \item $\deg(T_i\cdot T_i\cdot[F])=1$ for both $i=1,2$.
\end{enumerate}
In case (1), Proposition \ref{4p1} shows that $(T_1,T_2,[F])$ forms a regular fibration, giving the first alternative of Theorem \ref{Thm_1}.

It remains to show that case (2) yields the second alternative of Theorem \ref{Thm_1}.
For each integer $d\ge2$, set
\[
\mathrm{P}(d):=
\mathrm{conv}\bigl(
\{0\}\cup\{\pm e_i:1\le i\le d\}
\cup\{e_i-e_j:1\le i\neq j\le d\}
\bigr)
\]
and define
\[
\mathcal{M}(d):=
\left\{
(M_1,M_2)\in\mathrm{TR}(\R^d)\times\mathrm{TR}(\R^d)
\;\middle|\;
\begin{gathered}
\dim\bigl(\mathrm{Newt}(M_1)+\mathrm{Newt}(M_2)\bigr)=d,\\
\mathrm{Newt}(M_i)\subseteq\mathrm{P}(d)
\quad(i=1,2)
\end{gathered}
\right\}.
\]
The set $\mathcal{M}(d)$ is finite: $\mathrm{P}(d)$ contains only
finitely many lattice points, and each tropical regular function
is determined by its Newton polytope.

For each $(M_1,M_2)\in\mathcal{M}(d)$, let
$\mathcal{T}_{M_1,M_2}(d)$ be the set of tropical $2$-cycles $[G]$
in $\R^d$ satisfying
\[
\deg(M_i\cdot M_j\cdot[G])=1
\quad(i,j\in\{1,2\}).
\]
By Theorem \ref{4L4}, this set is finite. Hence the minimal-model class
\[
\mathcal{MD}_n:=
\bigcup_{d=2}^n
\bigcup_{(M_1,M_2)\in\mathcal{M}(d)}
\mathcal{T}_{M_1,M_2}(d)
\]
is a finite collection of tropical $2$-cycles.

Finally, Proposition \ref{4p2} gives a saturated $\Z$-projection
$\pi:\R^n\to\R^m$, with $2\le m\le n$, and a pair
$(M_1,M_2)\in\mathcal{M}(m)$ such that $\pi^*M_i=T_i$ for $i=1,2$.
Part (4) of that proposition gives
\[
\pi_*[F]\in\mathcal{T}_{M_1,M_2}(m)\subseteq\mathcal{MD}_n.
\]
This proves the second alternative and completes the proof.
\end{proof}

\begin{remark}
Although Theorem~\ref{Thm_1} does not fully answer
Question~\ref{Q}, it provides a powerful geometric engine:
for each regular tuple on a tropical regular $2$-cycle satisfying
the Hodge index inequality, it is either a regular fibration
or a saturated $\Z$-projection to one of finitely many
minimal models, for each fixed ambient dimension.
Understanding the relationship between these minimal models
and tropical smooth cycles is therefore a key step towards
answering Question~\ref{Q}.
In particular, how far can tropical minimal-model cycles
be from tropical smoothness, and can the minimal models
always be chosen to be tropically smooth?
These are still non-trivial open questions.
\end{remark}
\section{Regular Tropicalizations}\label{S5}

In this section, we aim to show that conjecture \hyperref[trop]{A.Esterov} holds for all regular tropicalizations $\Trop(U)$ with generic coefficients.
\begin{theorem}\label{Thm_2}
Let $U=V(f_1,\dots, f_{n-k})\subseteq (\C^*)^n$ be a $k$-dimensional very affine variety with generic coefficients. If $[F]:=[\Trop(U)]$ is regular, then any regular tuple $(T_1,\dots,T_k,[F])$ determines a saturated $\Z$-linear projection $\pi:\R^n\twoheadrightarrow \R^m$ such that one of the following holds:

\begin{enumerate}
    \item either the regular tuple $(T_1,\dots,T_k,[F])$ forms a regular fibration, and the base $\pi|F|$ is a support of a Bergman fan in $\R^m$,

    \item or $\pi_*[F]$ is a Bergman $k$-cycle in $\R^m$. There are tropical regular functions $M_1,\dots,M_k$ in $\mathrm{TR}(\R^m)$ such that
    $$T_i=\pi^*M_i\ \ (1\le i\le k),\qquad \deg(M_1\cdots M_k\cdot\pi_*[F])=1.$$
\end{enumerate} 
\end{theorem}

Let $U=V(f_1,\ldots,f_{n-k})\subseteq(\C^*)^n$, and put
$P_j:=\mathrm{Newt}(f_j)$. By \cite[Corollary~4.6.11]{MaclaganSturmfels15}, if the coefficients are sufficiently general for these fixed Newton polytopes, then
\begin{equation}
\label{rt:tropicalization}
    [\Trop(U)]=\Sigma_{P_1}\cap_{\mathrm{st}}\cdots\cap_{\mathrm{st}}\Sigma_{P_{n-k}}.
\end{equation}
 Here $\Sigma_{P_i}$ are hypersurfaces of $f_i$ (cf. Example \ref{ex:normal-fan}). Applying Theorem \ref{Bernstein} to the regularity of $[\Trop(U)]$, we obtain a tuple of lattice polytopes $Q_1,\dots ,Q_k\subseteq (\Z^n)^\vee$ satisfying
$$\mathrm{MV}(Q_1,\dots,Q_k,P_1,\dots,P_{n-k})=1.$$
Esterov-Gusev's classification theorem on mixed volume $1$ tuple of polytopes provides a powerful engine of our study:
\begin{theorem}[Esterov--Gusev \cite{EsterovGusev}]\label{E}
A collection $A$ of $n$ lattice polytopes in $(\mathbb R^n)^\vee$ has mixed volume one if and only if their mixed volume is nonzero and there exists $k>0$ such that up to translations, $k$ of the polytopes are faces of the same $k$-dimensional volume $1$ lattice simplex $\Delta_k$ in a $k$-dimensional rational subspace $U^\vee\subset(\mathbb R^n)^\vee$, and the images of the other $n-k$ polytopes under the projection $\mathbb (R^n)^\vee\to(\mathbb R^n)^\vee/U^\vee$ have mixed volume one.
\end{theorem}

Let $(P_1,\ldots,P_n)$ be a tuple of lattice polytopes in $(\R^n)^\vee$
with $\mathrm{MV}(P_1,\ldots,P_n)=1$.
Repeated application Theorem \ref{E}
 gives a partition of this tuple into
nonempty blocks $E_1,\ldots,E_s$. Write $n(i):=|E_i|$ for the number
of polytopes in $E_i$. Firstly, the polytopes in $E_1$ are faces of a unimodular simplex $\Delta_{n(1)}$ in an $n(1)$-dimensional rational subspace $U_1^\vee\subseteq(\R^n)^\vee$ up to translations while quotienting by $U_1^\vee$ gives a tuple of mixed
volume $1$ in $(\R^n)^\vee/U_1^\vee$.
Applying the theorem again, the images of the polytopes in $E_2$ are faces of a unimodular simplex $\Delta_{n(2)}$ in an
$n(2)$-dimensional rational subspace up to translations.
\[
    \overline{U_2^\vee}\subseteq(\R^n)^\vee/U_1^\vee.
\]
Let $U_2^\vee$ be the inverse image of $\overline{U_2^\vee}$ under the
quotient map. Then $U_1^\vee\subsetneq U_2^\vee$ and
$U_2^\vee/U_1^\vee=\overline{U_2^\vee}$, so
$\dim U_2^\vee=n(1)+n(2)$.
Continuing in this way gives the following partition and flag:

\begin{definition}[E-G partition]
\label{rt:EG_partition}
An \textit{E-G partition} of $(P_1,\ldots,P_n)$ consists of a partition
into nonempty blocks $E_1,\ldots,E_s$, with $n(i):=|E_i|$, together
with a flag of rational subspaces
\[
    0=U_0^\vee\subsetneq U_1^\vee\subsetneq\cdots
      \subsetneq U_s^\vee=(\R^n)^\vee,
    \qquad \dim(U_i^\vee/U_{i-1}^\vee)=n(i),
\]
such that, for every $i$, the images of the polytopes in $E_i$ under
$(\R^n)^\vee\to(\R^n)^\vee/U_{i-1}^\vee$ are faces of a unimodular simplex in $U_i^\vee/U_{i-1}^\vee$. Each $E_i$ is called an \textit{E-G block}.
\end{definition}

In particular, $\dim U_i^\vee=n(1)+\cdots+n(i)$ and
$n(1)+\cdots+n(s)=n$.
The correspondence between the blocks and the flag is shown below.

\begin{center}
\begin{tikzpicture}[>=stealth, every node/.style={inner sep=3pt}]
    \node (E1) at (0,1.1) {$E_1$};
    \node (E2) at (2.4,1.1) {$E_2$};
    \node at (4.8,1.1) {$\cdots$};
    \node (Es) at (7.2,1.1) {$E_s$};

    \node (U1) at (0,0) {$U_1^\vee$};
    \node at (1.2,0) {$\subsetneq$};
    \node (U2) at (2.4,0) {$U_2^\vee$};
    \node at (3.6,0) {$\subsetneq$};
    \node at (4.8,0) {$\cdots$};
    \node at (6.0,0) {$\subsetneq$};
    \node (Us) at (7.2,0) {$U_s^\vee$};
    \node[anchor=west] at (Us.east) {$=(\R^n)^\vee$};

    \draw[->] (E1.south) -- (U1.north);
    \draw[->] (E2.south) -- (U2.north);
    \draw[->] (Es.south) -- (Us.north);
\end{tikzpicture}
\end{center}

Now we fix a regular tuple $(T_1,\dots,T_k,[F])$ where $[F]=[\Trop(U)]=\Sigma_{P_1}\cap_{\mathrm{st}}\dots\cap_{\mathrm{st}}\Sigma_{P_{n-k}}$ for some lattice polytopes $P_i\subseteq (\R^n)^\vee$ and set $Q_i:=\mathrm{Newt}(T_i)$. The discussion above yields an \textit{E-G} partition for the tuple $(Q_1,\dots,Q_k,P_1,\dots,P_{n-k})$. Let $E_t$ be the \textbf{first} block containing a $Q_i$ and let $m=n(1)+\dots +n(t)$. After reordering the indices, suppose that $Q_1,\dots, Q_r$ are the $Q_i$'s in $E_t$, and the union $E_1\sqcup\dots\sqcup E_t$ consists of lattice polytopes
$$Q_1,\dots,Q_r,P_1,\dots,P_{m-r}.$$
Note that $r\le k$ and $m\le n$.

\begin{proposition}\label{5p1}
Assume $r=k$, and put $m:=\dim U_t^\vee$.
There is a saturated $\Z$-linear projection
$\pi:\R^n\twoheadrightarrow\R^m$ such that $\pi_*[F]$ is a Bergman
$k$-cycle in $\R^m$.
Moreover, there are tropical regular functions
$M_1,\ldots,M_k\in\mathrm{TR}(\R^m)$ such that
\[
    T_i=\pi^*M_i\quad(1\leq i\leq k),
    \qquad \deg(M_1\cdots M_k\cdot\pi_*[F])=1.
\]
\end{proposition}

\begin{proof}
Let $\overline{\Delta}_t\subseteq U_t^\vee/U_{t-1}^\vee$ be a
unimodular simplex such that the images of the polytopes in $E_t$ are faces of $\overline{\Delta}_t$. Then we
translate this simplex so that one vertex is the origin, lift its edge basis to $U_t^\vee$ and add a lattice basis of
$U_{t-1}^\vee$.
These vectors form a lattice basis of $U_t^\vee$.
Their convex hull with the origin is a unimodular $m$-simplex
$\Delta_t$ whose image in the quotient is $\overline{\Delta}_t$.

Choose the lattice coordinates given by this basis, so that
$\Delta_t=\mathrm{conv}(0,e_1,\ldots,e_m)$, and use the dual
coordinates on $U_t$.
Let
\[
    \pi:\R^n\twoheadrightarrow U_t\cong\R^m
\]
be dual to the inclusion $\iota:U_t^\vee\hookrightarrow(\R^n)^\vee$.
A lattice basis of $U_t^\vee$ extends to one of $(\Z^n)^\vee$, so
$\pi$ is saturated.
Since $r=k$, all $Q_i$ belong to $E_t$.
Thus $Q_i\subseteq U_t^\vee$.
Let $M_i$ be the support function of $Q_i$ on $U_t$.
Then $M_i\in\mathrm{TR}(U_t)$ and $T_i=\pi^*M_i$ for every $i$.

Set $H:=\max(0,x_1,\ldots,x_m)\in\mathrm{TR}(\R^m)$, so that
$\mathrm{Newt}(\pi^*H)=\iota(\Delta_t)$.
The key observation is
\[
    \mathrm{MV}\left(
        \underbrace{\iota(\Delta_t),\ldots,\iota(\Delta_t)}_{k\text{ times}},
        P_1,\ldots,P_{n-k}\right)=1.
\]
Indeed, replace every $Q_i$ by $\iota(\Delta_t)$ and keep the partition
$E_1,\dots,E_s$ and its flag. All blocks except $E_t$ are unchanged.
In $U_t^\vee/U_{t-1}^\vee$, the images of the $P_i$'s in $E_t$
remain unchanged, while each replacement has image $\overline{\Delta}_t$.
Hence the images in the modified block $E_t$ are still
faces of $\overline{\Delta}_t$ up to translations. Thus the same partition and flag satisfy the criterion
in Theorem \ref{E}, so the modified tuple also has mixed volume $1$.

Theorem \ref{Bernstein} and the projection formula now give
\[
    \deg\bigl(H^k\cdot\pi_*[F]\bigr)
      =\deg\bigl((\pi^*H)^k\cdot[F]\bigr)=1.
\]
Hence $\pi_*[F]$ is an effective $k$-cycle of degree $1$.
Theorem \ref{AlexFink} shows that it is a Bergman
$k$-cycle.
Finally, the projection formula also gives
\[
    \deg(M_1\cdots M_k\cdot\pi_*[F])
      =\deg(T_1\cdots T_k\cdot[F])=1
\]
as required.\\
\end{proof}

\begin{lemma}\label{5L1}
  Let $U^\vee\subset(\R^n)^\vee$ be an $m$-dimensional rational subspace with inclusion
  $\iota:U^\vee\hookrightarrow(\R^n)^\vee$. Let
  $\overline{P}_1,\dots,\overline{P}_k$ $(k\le m)$ be lattice polytopes in $U^\vee$ and set
  $P_i:=\iota(\overline{P}_i)\subset(\R^n)^\vee$.  Define
  \[
    F \;:=\; \Sigma_{P_1}\cap_{\mathrm{st}}\cdots\cap_{\mathrm{st}}\Sigma_{P_k},
    \qquad
    \overline{F} \;:=\; \Sigma_{\overline{P}_1}\cap_{\mathrm{st}}\cdots
        \cap_{\mathrm{st}}\Sigma_{\overline{P}_k},
    \qquad
    U^\perp:=\{x\in\R^n\mid \langle x,u\rangle=0\ \forall\, u\in U^\vee\}.
  \]
  Then $F=\overline{F}\times \mathcal{U}^\perp$ where $\mathcal{U}^\perp$ is a complete fan supported on $U^\perp$.
\end{lemma}

\begin{proof}
After choosing coordinates $(x_1,\dots,x_m,y_1,\dots,y_{n-m})$ on $\R^n$ so that
  $U^\vee=V(\langle y_1,\dots,y_{n-m}\rangle)$. For each $i$, let
  $\overline{f}_i\in \C[x_1^{\pm},\dots,x_m^{\pm}]$ be a Laurent polynomial with
  generic coefficients and $\mathrm{Newt}(\overline{f}_i)=\overline{P}_i$, and
  let $f_i\in \C[x_1^{\pm},\dots,x_m^{\pm},y_1^{\pm},\dots,y_{n-m}^{\pm}]$
  denote the same polynomial viewed in the larger Laurent ring, so that
  $\mathrm{Newt}(f_i)=P_i$.  Then we have
  \[
    F \;=\; \mathrm{Trop}\bigl(V(f_1,\dots,f_k)\bigr).
  \]
  Since $f_i$ does not involve the variables $y_1,\dots,y_{n-m}$,
  \[
    V(f_1,\dots,f_k)
      \;=\; V(\overline{f}_1,\dots,\overline{f}_k)\times \mathbb{T}^{n-m}_y.
  \]
  As tropicalization respects products,
  \[
    F \;=\; \mathrm{Trop}\bigl(V(\overline{f}_1,\dots,\overline{f}_k)\bigr)
            \times \mathrm{Trop}(\mathbb{T}^{n-m}_y)
      \;=\; \overline{F}\times \mathcal{U}^\perp,
  \]
  where $\mathrm{Trop}(\mathbb{T}^{n-m}_y)=\colon \mathcal{U}^\perp$ as required.\\
\end{proof}

\begin{proposition}\label{5p2}
    Assume $0<r<k$, and put $m:=\mathrm{dim}U_t^\vee$. There is a saturated $\Z$-linear projection $\pi:\R^n\twoheadrightarrow \R^m$ such that $\pi|F|$ is a support of a Bergman fan in $\R^m$ and $(T_1,\dots,T_k,[F])$ forms a regular fibration.
\end{proposition}

\begin{proof}
We admit the same saturated $\Z$-linear projection $\pi:\R^n\twoheadrightarrow \R^m$ in the proof of Proposition \ref{5p1} above. Recall that we may reorder the polytopes associated with $U_t^\vee$ as follows
$$\underbrace{Q_1,\dots Q_r}_{T_i\in\mathrm{TR}(\R^m)},\qquad \underbrace{P_1,\dots,P_{m-r}}_{[F]=\Sigma_{P_1}\cap_{\mathrm{st}}\dots\cap_{\mathrm{st}}\Sigma_{P_{n-k}}}$$
Now we define $\overline{P_i}$ (resp. $\overline{Q_j})$ to be the lattice polytope in $U_t^\vee$ such that $\iota(\overline{P_i)}=P_i$ for all $1\le i\le m-r$ (resp. $\iota(\overline{Q_j})=Q_j$ for all $1\le j\le r$). The Esterov-Gusev decomposition yields that their mixed volume is one:
$$\mathrm{MV}_{U_t^\vee}(\overline{Q_1},\dots,\overline{Q_r},\overline{P_1},\dots,\overline{P_{m-r}})=1.$$
The same argument in the proof of Proposition \ref{5p1} allows us to replace $\overline{Q_i}$ by $\overline{\Delta}_t$ and implies that the stable intersection $$[G]:=\Sigma_{\overline{P_1}}\cap_{\mathrm{st}}\dots\cap_{\mathrm{st}}\Sigma_{\overline{P_{m-r}}}$$
is an $r$-dimensional Bergman cycle in $U_t$.
Now we define two new tropical cycles in $\R^n$:
$$X=\Sigma_{P_1}\cap_{\mathrm{st}}\dots\cap_{\mathrm{st}}\Sigma_{P_{m-r}},\qquad Y=T_{r+1}\cdots T_{k}\cdot[F]=X\cap_{\mathrm{st}}\underbrace{\Sigma_{Q_{r+1}}\cap_{\mathrm{st}}\dots\cap_{\mathrm{st}}\Sigma_{P_{n-k}}}_{\text{polytopes not in $E_1\sqcup\dots\sqcup E_t$}}.$$
We aim to show that $\pi|F|=|G|$, so that the image is a support of an $r$-dimensional Bergman fan in $U_t\cong \R^m$ and $(T_1,\dots,T_k,[F])$ forms a regular fibration immediately.

Notice that the supports of $X,[F],Y$ have the following relation:
$$|Y|\subseteq |F|\subseteq |X|$$
since $[F]$ can be viewed by taking tropical intersections on $X$, and $Y$ can be viewed by taking tropical intersections on $[F]$. Then the relation above follows from the definition of tropical intersection product (or the stable intersection) directly. As a consequence, the supports of their images under $\pi$ have
$$\pi|Y|\subseteq \pi|F|\subseteq \pi|X|$$
so it is sufficient to prove that $\pi|X|=|G|=\pi|Y|$. We first denote $D_i$ to be the tropical regular function on $\R^n$ associated with $P_i$ for all $1\le i\le m-r$. In other words, $D_i$ is the support function of ${P_i}$. We may define the analogue tropical regular function $\overline{D_i}$ on $U_t$ which is associated with $\overline{P_i}$ for all $1\le i\le m-r$. Then $P_i=\iota(\overline{P_i})$ implies that $\pi^*\overline{D_i}=D_i$ for all $1\le i\le m-r$. Now we compute them separately: 

On the one hand, Lemma \ref{5L1} tells us that $X=G\times\mathcal{U}^\perp$ for some complete tropical cycle $\mathcal{U}^\perp$ supported in $U_t^\perp$. Therefore $\pi|X|=\pi|G\times \mathcal{U}^\perp|=|G|$ follows from $\ker\pi=U_t^\perp$ immediately.

On the other hand, we denote $Y=X\cap_{\mathrm{st}}\mathcal{Z}$ where $\mathcal{Z}$ is the stable intersection of hypersurfaces whose Newton polytopes are not in $E_1\sqcup\dots\sqcup E_t$. Note that $\mathcal{Z}$ is a non-zero $m$-dimensional tropical cycle. Then the tropical projection formula yields
\begin{align*}
    \pi_*Y&=\pi_*(X\cap_{\mathrm{st}}\mathcal{Z})\\
    &=\pi_*(D_1\cdots D_{m-r}\cdot\mathcal{Z})\\
    &=\pi_*(\pi^*\overline{D_1}\cdots\pi^*\overline{D_{m-r}}\cdot\mathcal{Z})\\
    &=\overline{D_1}\cdots\overline{D_{m-r}}\cdot\pi_*\mathcal{Z}.
\end{align*}
As a consequence, $\pi_*\mathcal{Z}$ is a $m$-dimensional tropical cycle supported on $U_t\cong\R^m$. Therefore $|\pi_*Y|=|G|$. In other words, $\pi|Y|=|G|$.

Combine these together, we conclude that $|G|=\pi|Y|$ which completes the proof.\\
\end{proof}

\begin{proof}[proof of Theorem \ref{Thm_2}]
   For any regular tuple $(T_1,\dots,T_k,[F])$ where $$[F]=\Sigma_{P_1}\cap_{\mathrm{st}}\dots\cap_{\mathrm{st}}\Sigma_{P_{n-k}}$$ we choose an \textit{E-G} decomposition of $(Q_1,\dots,Q_k,P_1,\dots,P_{n-k})$. Then the proof follows directly from Propositions \ref{5p1} and \ref{5p2}.
\end{proof}

\begin{remark}
Theorem~\ref{Thm_2} answers Question~\ref{Q} for the regular
tropicalizations considered in this section by relating them
to tropical smooth cycles.
For every regular tuple $(T_1,\ldots,T_k,[F])$, there is a
saturated $\mathbb{Z}$-projection $\pi$ whose
set-theoretic image $\pi(|F|)$ supports a Bergman cycle $[G]$.
In case~(1), $[G]$ is the base of a regular fibration, and
some of the functions $T_i$ are pulled back from tropical
regular functions on its ambient space.
In case~(2), $[G]=\pi_*[F]$, and all the functions $T_i$
can be obtained in this way.

A stronger goal is to describe these projections intrinsically
in terms of the combinatorial geometry of $[F]$, without
first specifying a regular tuple.
In fact, this holds for regular $1$-cycles, see \cite{Li25}.
Obtaining an analogous description in higher dimensions
remains a nontrivial open problem.
\end{remark}
\section{Tropical cycles with an unmixed sequence}\label{S6}

This section concerns the following question:
\begin{center}Let $[F]$ be a tropical $k$-cycle in $\mathbb{R}^{n}$, and let $T$ be a tropical regular function. Can one give a combinatorial classification of the tuples $(T,\dots,T,[F])$ arising from an unmixed sequence according to $\deg(\underbrace{T\cdots T}_{k\text{ times}}\cdot[F])$ equals $0$ or $1$?
\end{center}

The degree-zero case for unmixed sequences is a natural question in its own right: before studying intersection number one, it is natural to understand the case of intersection number zero.

Our interest in the unmixed degree-one case has two motivations. First, we aim to establish the conjecture \hyperref[trop]{A.~Esterov} holds for unmixed regular tuples. Second, Theorem~\ref{AlexFink} characterizes Bergman fans by an unmixed intersection number of one with respect to the standard tropical regular function $\max(0,x_1,\ldots,x_n)$. Therefore, a natural first step towards extending this characterization of tropical smoothness is to replace the standard function by an arbitrary tropical regular function $T$ with $L(T)=\{0\}$.

We now recall the basic material on matroid simplification that will be needed later.
\begin{definition}[Parallel elements,\cite{Oxley},\cite{MaclaganSturmfels15}]
Let $M=(E,\rho)$ be a matroid. Two distinct non-loop elements $e,f\in E$ are called parallel if $\{e,f\}$ is a circuit of $M$. Equivalently, $e$ and $f$ are parallel if $\rho(\{e,f\})=1$.

The relation generated by $e\sim f$ if either $e=f$ or $e$ and $f$ are parallel is an equivalence relation on the set of non-loop elements of $M$. Its equivalence classes are called the parallel classes of $M$.
\end{definition}

\begin{definition}[Matroid simplification,\cite{Oxley},\cite{MaclaganSturmfels15}]
Let $M$ be a matroid on a finite ground set $E$. The simplification of $M$, denoted by $M^\mathrm{si}$, is obtained by deleting all loops of $M$ and then replacing each parallel class by a single element.

Equivalently, let $\overline{E}$ be the set of parallel classes of non-loop elements of $M$. For a subset $A\subseteq \overline{E}$, define
$\widetilde A:=\bigcup_{P\in A}P\subseteq E$. The rank function of $M^\mathrm{si}$ is given by
$\rho_{M^\mathrm{si}}(A):=\rho_M(\widetilde A)$.
Thus $M^\mathrm{si}$ is a loopless matroid with no parallel elements, and it has the same rank as $M$ after removing loops.
\end{definition}

\begin{proposition}\label{simplification} Let $M$ be a loopless matroid on a finite ground set $E$ with set of parallel classes $\mathcal P(M)$, and let $M^{\mathrm{si}}$ be its simplification, with ground set $\mathcal P(M)$. The $\mathbb Z$-linear surjection $$\pi\colon\mathbb R^E/\mathbb R\mathbf 1\twoheadrightarrow\mathbb R^{\mathcal P(M)}/\mathbb R\mathbf 1,\qquad (w_e)_{e\in E}\longmapsto(w_{e_P})_{P\in\mathcal P(M)}$$ (for any choice of representatives $e_P\in P$) satisfies $\pi_*\bigl(\mathrm{Berg}(M)\bigr)=\mathrm{Berg}(M^{\mathrm{si}})$. In particular, if $\mathrm{Berg}(M)$ has reduced ambient dimension $m$, then $|\mathcal P(M)|=m+1$. \end{proposition}

\begin{proof} If $e,f$ are parallel, then $\{e,f\}$ is a circuit, so every $w\in\mathrm{Berg}(M)$ has $w_e=w_f$. Thus $\mathrm{Berg}(M)$ lies in the subspace $V_0\subseteq\mathbb R^E/\mathbb R\mathbf 1$ and $\pi$ is independent of the chosen representatives.

Since $M$ is loopless, every flat $F$ of $M$ is a union of parallel classes (if $e\in F$ and $f$ is parallel to $e$, then $f\in\mathrm{cl}(F)=F$). Writing $F^{\mathrm{si}}:=\{P\in\mathcal{P}(M)\mid P\subseteq F\}$, the maps $F\mapsto F^{\mathrm{si}}$ and $A\mapsto\widetilde A:=\bigcup_{P\in A}P$ are mutually inverse, order-preserving lattice isomorphisms between the flats of $M$ and those of $M^{\mathrm{si}}$. Note that $F^{\mathrm{si}}$ is a flat of $M^{\mathrm{si}}$ since $\rho_{M^{\mathrm{si}}}(F^{\mathrm{si}}\cup\{P\})=\rho_M(F\cup P)>\rho_M(F)$ for $P\not\subseteq F$.

The cones of $\mathrm{Berg}(M)$ are generated by chains of proper flats, with rays $-\mathbf 1_F$. Notice that $e_P\in F\iff P\in F^{\mathrm{si}}$, we have $\pi(-\mathbf 1_F)=\mathbf -1_{F^{\mathrm{si}}}$, and the flat isomorphism sends chains to chains; hence $\pi$ maps the cones of $\mathrm{Berg}(M)$ bijectively onto those of $\mathrm{Berg}(M^{\mathrm{si}})$. On $V_0$ the map sends $-\mathbf 1_P\mapsto\mathbf e_P$, so it restricts to a lattice isomorphism $V_0\cap(\mathbb Z^E/\mathbb Z\mathbf 1)\xrightarrow{\sim}\mathbb Z^{\mathcal P(M)}/\mathbb Z\mathbf 1$; thus every pushforward multiplicity is $1$, and since all weights equal $1$, $\pi_*(\mathrm{Berg}(M))=\mathrm{Berg}(M^{\mathrm{si}})$. Finally $\mathbb R^{\mathcal P(M)}/\mathbb R\mathbf 1\cong\mathbb R^{|\mathcal P(M)|-1}$ has dimension $m$, so $|\mathcal P(M)|=m+1$. \end{proof}

\subsection{$\deg(T^k\cdot[F])=0$}

The principal result of this subsection is Theorem~\ref{TF=0_thm}, which provides a necessary and sufficient condition for the tuple $(T,\dots,T,[F])$ to have degree zero.

\begin{lemma}\label{TF=0_lemma1}
    Let $[F]$ be a non-zero tropical $k$-cycle in $\R^n$ and let $T$ be a tropical regular function. If $T\cdot [F]$ is a zero tropical $k-1$ cycle, then $F\subset L(T)$ where $L(T)$ is the lineality space of $T$ (cf. Definition \ref{TR}).
\end{lemma}

\begin{proof} Write $T=\max(0,l_1,\dots,l_m)$. Let $\tau\in F^{(k-1)}$, and for each maximal cone $\sigma>\tau$ let $u_{\sigma/\tau}$ be the primitive generator of $\sigma$ modulo $\langle\tau\rangle$. Now we choose $F$ to be a representative of $[F]$ such that $T$ is conewise-linear over $F$. This suggests a map $ F\to2^{\{0,l_1,\dots,l_m\}}$ defined as follows: for any cone (not necessarily maximal) $\zeta\in F$,
$$\zeta\mapsto S_\zeta\subseteq\{0,l_1,\dots,l_m\}$$
such that $T_\zeta(v)=l_i(v)$ for all $l_i\in S_\zeta$ and $v\in\zeta$. The definition of the tropical intersection product gives $$\omega_{T\cdot[F]}(\tau)=\sum_{\sigma>\tau}T_\sigma\bigl(\omega_F(\sigma)\,u_{\sigma/\tau}\bigr)-T_\tau(\sum_{\sigma>\tau}\omega_F(\sigma)\,u_{\sigma/\tau}).$$
We may replace the restricted function $T_\tau$ by an arbitrary linear functional $l_\tau\in S_\tau$. Hence we obtain
$$\omega_{T\cdot[F]}(\tau)=\sum_{\sigma>\tau}\bigl(T_\sigma\bigl(\omega_F(\sigma)\,u_{\sigma/\tau}\bigr)-l_\tau(\omega_F(\sigma)\,u_{\sigma/\tau}) \bigr).$$
On the one hand, $T$ is taking maximum among $\{0,l_1,\dots,l_m\}$, so $$T_\sigma(\omega_F(\sigma)u_{\sigma/\tau})\ge l_\tau(\omega_F(\sigma)u_{\sigma/\tau})$$for all $\sigma>\tau$ since $\omega_F(\sigma)>0$. On the other hand, $\omega_{T\cdot[F]}(\tau)=0$ forces
$$T_\sigma(\omega_F(\sigma)u_{\sigma/\tau})= l_\tau(\omega_F(\sigma)u_{\sigma/\tau})$$
for all $\sigma>\tau$.

It is sufficient to prove the assertion for codimension-one connected
tropical cycles. Indeed, let \(C_1,\ldots,C_r\) be the codimension-one
connected components of \(F\). Then
\[
[F]=\sum_{i=1}^r [C_i].
\]
By the \(\mathbb Z\)-linearity of the tropical divisor product in the cycle
variable, we have
\[
T\cdot [F]=\sum_{i=1}^r T\cdot [C_i].
\]
Moreover, each \(T\cdot [C_i]\) is an effective tropical cycle, since \(T\) is
a tropical regular function. Hence \(T\cdot [F]=0\) if and only if
\(T\cdot [C_i]=0\) for every \(i\). Therefore it remains to prove the claim
under the additional assumption that \(F\) is codimension-one connected.

Assume now that \(F\) is codimension-one connected. The argument above
shows that \(T_\sigma=T_{\sigma'}\) whenever two maximal cones
\(\sigma,\sigma'\) meet along a codimension-one cone. Hence, by
codimension-one connectedness, there exists a single integral linear form
\(l\) such that
\[
T|_\sigma=l
\]
for every maximal cone \(\sigma\) of \(F\). Thus \(T=l\) on \(|F|\).

Since $F$ is a balanced codimension-one connected fan with strictly positive weights, its support is not contained in any proper half-space of its real span. In other words,
$$\mathrm{cone}(|F|)=\mathrm{span}_\R(|F|)=\colon V,$$
which is a rational linear subspace of $\R^n$. If $l\neq 0$ on $V$, then $l$ takes a negative value at some point of $V$. Equivalently
$$0>l(v)=l(\sum a_i\rho_i)=\sum a_i T(\rho_i)$$
where $v=\sum a_i\rho_i$ for some $\rho_i\in|F|$ and $a_i>0$. However, this contradicts the non-negativity of $T=\max(0,l_1,\dots,l_m)$. Hence $l|_V=0$.

It follows that $T=0$ on $|F|$, i.e. $|F|\subset L_0$ where $L_0$ denotes the cone $\{x\in\R^n\mid T(x)=0\}$. Note that any linear subspace of $L_0$ is contained in $L(T)$. Therefore $F\subset L_0$ implies that $V\subseteq L_0$ whence $V\subseteq L(T)$. So one concludes that $F\subset L(T)$ as required.

\end{proof} 

\begin{lemma}\label{TF=0_lemma_2}
    Let $[F]$ be a non-zero tropical $k$-cycle in $\R^n$ and let $T:=\max(0,l_1,\dots,l_m)$ be a tropical regular function whose lineality space $L(T)=\{0\}$. Then 
    $$\deg(T^k\cdot[F])>0.$$
\end{lemma}

\begin{proof}
 We will prove it by induction on the dimension. First, let $[C]$ be a non-zero tropical $1$-cycle in $\R^n$. The balancing condition ensures that there exists at least one ray of $C$ which has a primitive generator $v$ with $T(v)>0$, otherwise $C$ is contained in a convex cone of $\R^n$ which contradicts the balancing condition. Hence $\deg(T\cdot[C])\ge T(v)>0$.

Now fix an integer $h>1$ and assume the statement holds for every non-zero tropical cycle in $\R^n$ of dimension strictly less than $h$. Let $[G]$ be a non-zero tropical $h$-cycle in $\R^n$, and let $T$ be a tropical regular function with lineality space $L(T)=\{0\}$. By Lemma \ref{TF=0_lemma1}, the condition $L(T)=\{0\}$ implies that $T\cdot[G]$ is a non-zero tropical $(h-1)$-cycle. The induction hypothesis then gives
$$
\deg\bigl(T^{h-1}\cdot(T\cdot[G])\bigr)>0,
$$
which is the desired conclusion. This completes the induction, and therefore $\deg(T^k\cdot[F])>0$, as required.
\end{proof}

\smallskip

Now we introduce one of the main results of this section:

\begin{theorem}\label{TF=0_thm}
Let $[F]$ be a tropical $k$-fan in $\mathbb R^n$ and let $T$ be a tropical regular function. Then $\deg(T^k\cdot [F])=0$ if and only if
$L(T)\cap \langle\sigma\rangle\neq 0$ for every maximal cone $\sigma\in F^{(k)}$, where
$\langle\sigma\rangle:=\operatorname{span}_{\mathbb R}(\sigma)$.
\end{theorem}

\begin{proof}
Set $L:=L(T)$, and let $\pi:\mathbb R^n\to \mathbb R^n/L$ be the quotient map. Since $T$ is constant along $L$, there exists a tropical regular function $T'$ on $\mathbb R^n/L$ such that $T=\pi^*(T')$. By construction, $L(T')=\{0\}$.

By the projection formula, we have $\pi_*(T^k\cdot [F])=T'^k\cdot \pi_*[F]$. Hence
$\deg(T^k\cdot [F])=0$ if and only if $\deg(T'^k\cdot \pi_*[F])=0$.

By Lemma \ref{TF=0_lemma_2}, since $L(T')=\{0\}$, the equality $\deg(T'^k\cdot \pi_*[F])=0$ holds if and only if $\pi_*[F]=0$ as a $k$-cycle. For an effective fan, this is equivalent to saying that every maximal cone of $F$ drops dimension under $\pi$. In other words, for every $\sigma\in F^{(k)}$, one has $\dim \pi(\sigma)<k$. Therefore we have
$L\cap \langle\sigma\rangle\neq 0$ for every maximal cone $\sigma$ of $F$.

Thus $\deg(T^k\cdot [F])=0$ if and only if $L(T)\cap\langle\sigma\rangle\neq 0$ for all maximal cones $\sigma\in F^{(k)}$, as claimed.\\
\end{proof}

\smallskip

\subsection{$\deg(T^k\cdot[F])=1$} In this subsection, we will show that the conjecture \hyperref[trop]{A.Esterov} holds for unmixed regular tuples $(\underbrace{T,\dots,T}_{\text{k times}},[F])$ where the lineality space $L(T)=\{0\}$.

\begin{theorem}\label{Thm_3}
    Let $[F]$ be a tropical $k$-cycle over $\R^n$ which spans $\R^n$, and let $T$ be a
    tropical regular function satisfying $\deg(T^k \cdot [F]) = 1$ and the lineality space $L(T) = \{0\}$ (cf. Definition \ref{TR}). Then
    there exists a $\mathbb{Z}$-isomorphism $\psi$ such that
    \begin{enumerate}
        \item $\psi_*[F]$ is a Bergman $k$-cycle spanning $\mathbb{R}^n$, and
        \item there exists a tropical regular function $M$ such that $\deg(M^k \cdot \psi_*[F]) = 1$ and
              $T = \psi^*(M)$.
    \end{enumerate}
    Moreover, $\mathrm{Newt}(T)$ is an $n$-simplex.
\end{theorem}

\begin{proof}
    Write the regular function as $T = \max(0, l_1, \dots, l_q)$, and define
    the $\mathbb{Z}$-linear map $\pi \colon \mathbb{Z}^n \to \mathbb{Z}^q$
    whose matrix representation is
    \[
        \pi \;=\;
        \begin{bmatrix}
            l_1 \\ \vdots \\ l_q
        \end{bmatrix}.
    \]
    Denoting the standard basis of $\mathbb{Z}^q$ by $e_1, \dots, e_q$, it follows that
    $\pi^*(\max(0, e_1, \dots, e_q)) = T$. Since $L(T) = \{0\}$ we have
    $\ker(\pi) = 0$; in particular, $q\geq n$. It follows that $\pi_*(F)$ is
    a non-zero $k$-fan over $\mathbb{R}^q$. Moreover, since $[F]$ spans
    $\mathbb{R}^n$, the fan $\pi_*(F)$ is contained in an $n$-dimensional
    subspace of $\mathbb{R}^q$.

    Applying the projection formula, we obtain
    \begin{align*}
        1 \;=\; \deg(T^k \cdot [F])
        &\;=\; \deg(\pi_*\bigl(\pi^*(\max(0,e_1,\dots,e_q))^k
               \cdot [F])) \\
        &\;=\; \deg(\max(0,e_1,\dots,e_q)^k
               \cdot \pi_*[F]).
    \end{align*}
    We conclude that $\pi_*[F]$ is a Bergman $k$-cycle over $\mathbb{R}^q$
    contained in a rational subspace of dimension $n$ via Theorem \ref{AlexFink}.

    After simplifying the matroid underlying $\pi_*(F)$, Proposition \ref{simplification} yields a $\Z$-linear surjection $\rho$ and a rank-$k+1$
    matroid $(\mathbf{M}, E)$ whose
    Bergman cycle $[\mathrm{Berg}(\textbf{M})]$ is isomorphic to $\pi_*[F]$ under $\rho$. Consider the $\Z$-linear map $\psi\colon=\rho\circ\pi\colon\R^n\to \R^n$, it follows that
    $$\psi_*[F]=[\mathrm{Berg}(\textbf{M})].$$
    We will show that $\psi$ is a $\Z$-isomorphism. It is enough to show that $\psi(\Z^n)=\Z^n$. We may assume that $\mathcal{F}$ is a representative of the tropical cycle $[F]$ such that $\psi_*(\mathcal{F})=\mathrm{Berg}(\textbf{M})$.

By the definition of pushforward of fans (cf.\ Definition~\ref{def:pushforward}),
the weight of $\tilde{\sigma}\in\psi_*(\mathcal{F})$ at a maximal cone $\sigma$ equals
$$\sum_{\tilde{\sigma}=\psi(\sigma)}
[\Z^n_{\psi(\sigma)}:\psi(\Z^n_{\sigma})]\cdot w_{\mathcal{F}}(\sigma).
$$
Since the weight function of a Bergman fan is identically $1$, one deduces
$$
[\Z^n_{\psi(\sigma)}:\psi(\Z^n_{\sigma})]=1
$$
for every maximal cone $\sigma$ of $\mathcal{F}$.

Now let $\R_{\ge 0}\cdot v$ be a ray of $\psi_*(\mathcal{F})$ with primitive
generator $v\in \Z^n$. Choose a maximal cone $\tau'$ of $\psi_*(\mathcal{F})$ containing
$\R_{\ge 0}\cdot v$, and a maximal cone $\tau$ of $\mathcal{F}$ with $\psi(\tau)=\tau'$. Then
$$
v\;\in\;\langle\tau'\rangle\cap \Z^n\;=\;\Z^n_{\tau'}\;=\;\Z^n_{\psi(\tau)}\;=\;\psi(\Z^n_{\tau}),
$$
so there exists $w\in\Z^n_{\tau}\subset\Z^n$ with $\psi(w)=v$. Since the primitive ray generators of a Bergman fan span its ambient lattice $\Z^n$, we conclude that $\Z^n\subseteq\psi(\Z^n)$, hence $\psi(\Z^n)\cong\Z^n$ as required.
\smallskip

To verify condition~(2), set $\phi := \psi^{-1}$ and apply the projection formula again
    \begin{equation*}
    \begin{aligned}
        \deg(T^k\cdot [F])
        \;&=\; \deg(T^k\cdot \phi_*[\mathrm{Berg}(\mathbf{M})])\\
        &=\; \deg(\phi_*(\phi^*(T)^k
              \cdot [\mathrm{Berg}(\mathbf{M})]))\\
        &=\; \deg(\phi^*(T)^k
              \cdot [\mathrm{Berg}(\mathbf{M})])\\
        &=\; 1.
        \end{aligned}
    \end{equation*}
    Hence $\psi_*[F] = [\mathrm{Berg}(\mathbf{M})]$ is a Bergman $k$-cycle
    spanning $\mathbb{R}^n$, and $M := \phi^*(T)$ is a regular sequence for
    $\mathrm{Berg}(\mathbf{M})$ satisfying $T = \psi^*(M)$ as required.
\smallskip

To see the remaining statement holds, we distinguish two cases.

    \medskip\noindent\textit{Case~1: $T^{k-1} \cdot [F]$ spans $\mathbb{R}^n$ as a $1$-curve.} The classification theorem for $1$-cycles established in \cite{Li25}
    gives that $T$ is a pullback of a tropical regular function $M$ of a Bergman $1$-fan $\mathcal{C}$ such that $\deg(M\cdot[\mathcal{C}])=1$ over $\mathbb{R}^n$. Therefore, $\mathrm{Newt}(T)$ is a simplex such that $\mathrm{dim}(\mathrm{Newt}(T))\le n$. Then $\mathrm{dim}(\mathrm{Newt}(M))=n$ follows directly from $L(T)=\{0\}$. Hence $\mathrm{Newt}(T)$ is an $n$-simplex as required.

    \medskip\noindent\textit{Case~2: $T^{k-1} \cdot [F]$ does not span $\mathbb{R}^n$.}
    Then there exists a functional $l \in (\mathbb{Z}^n)^\vee$ such that
    $T^{k-1} \cdot [F] \subset \ker(l)$, whence $\deg(\max(0,l) \cdot T^{k-1} \cdot [F]) = 0$.
    So
   \begin{equation}
        \deg(\max(0,l) \cdot T^{k-1} \cdot [F])
        \;=\;\deg( T^{k-1} \cdot \max(0,l) \cdot [F])
        \;=\; 0.
   \end{equation}
    Recall that $[F]$ is a $k$-cycle spanning $\mathbb{R}^n$, Lemma \ref{TF=0_lemma1} yields that the product $\max(0,l) \cdot [F]$ is a nonzero $k-1$-cycle over $\mathbb{R}^n$. The condition $L(T)=\{0\}$ allows us to apply Lemma \ref{TF=0_lemma_2} on the tuple $(\underbrace{T,\dots,T}_{k-1\text{ times}},\max(0,l)\cdot[F])$:
    $$\deg(\underbrace{T,\cdots,T}_{k-1\text{ times}}\cdot \max(0,l)\cdot[F])>0$$
 which contradicts (6.1).

    \medskip
    Combining both cases, we conclude that $\mathrm{Newt}(T)$ is an $n$-simplex as required.
\end{proof}

\begin{remark} The spanning hypothesis $\mathrm{span}_{\mathbb{R}}(|F|)=\mathbb{R}^{n}$ in Theorem \ref{Thm_3} cannot be omitted. Consider the tropical $1$-cycle $[F]$ in $\mathbb{R}^{3}$ whose rays are $$\rho_{1}=\mathbb{R}_{+}(1,0,0),\qquad\rho_{2}=\mathbb{R}_{+}(-1,1,0),\qquad\rho_{3}=\mathbb{R}_{+}(0,-1,0),$$ each with weight $1$. Since $(1,0,0)+(-1,1,0)+(0,-1,0)=(0,0,0)$, the balancing condition of Definition~\ref{tropical_fans} holds, so $[F]$ is a tropical $1$-cycle over $\mathbb{R}^{3}$; note that $\mathrm{span}_{\mathbb{R}}(|F|)$ is the plane $\{z=0\}$, so $[F]$ does not span $\mathbb{R}^{3}$. Let $$T=\max\bigl(0,\;x,\;x+y,\;x+\lambda z,\;\mu z\bigr),\qquad\lambda,\mu\in\mathbb{Z}_{>0}.$$ The four linear forms $0,\,x,\,x+y,\,\mu z$ already have only the origin as common zero, so $L(T)=\{0\}$. On the other hand, we conclude that $T|_{\{z=0\}}=\max(0,x,x+y)$; a direct computation of the stable intersection $V(T)\cap_{\mathrm{st}}F$ then gives $$\deg(T\cdot[F])=1,$$ so $(T,[F])$ is a regular tuple. However, $\mathrm{Newt}(T)=\mathrm{conv}\{0,(1,0,0),(1,1,0),(1,0,\lambda),(0,0,\mu)\}$ is a full-dimensional polytope with five vertices; it is therefore not a simplex, nor is $\mathbb{Z}$-isomorphic to one. Thus the conclusion of Theorem \ref{Thm_3} fails once the spanning hypothesis is dropped. \end{remark}

\begin{remark}
Theorem \ref{Thm_3} answers Question~\ref{Q} for unmixed
regular tuples with trivial lineality.
Indeed, if a tropical $k$-cycle $[F]$ spans $\mathbb{R}^n$
and admits a tropical regular function $T$ satisfying
\[
\deg(T^k\cdot[F])=1,
\qquad L(T)=\{0\},
\]
then a lattice isomorphism identifies $[F]$ with a Bergman
cycle. Thus $[F]$ itself is tropically smooth.
\end{remark}

\section*{Tool and computational resource disclosure}
In accordance with the Leiden Declaration on Artificial Intelligence and Mathematics, which the author has signed (3 June 2026), the author declares the use of automated tools in this paper as follows.

The main mathematical content and logical reasoning in this paper were completed independently by the author under the supervision of \textsc{Alexander Esterov}.

ChatGPT was used to simplify author's original proof of Theorem \ref{Thm_2}, optimize the language (including grammar and words), and check the proofs in the paper.


\begin{thebibliography}{99}
\bibitem{AHK18}
Karim~Adiprasito, June~Huh, and Eric~Katz,
\emph{Hodge theory for combinatorial geometries},
Ann. of Math. (2) \textbf{188} (2018), no.~2, 381--452.
\url{https://doi.org/10.4007/annals.2018.188.2.1}.

\bibitem{AksnesAminiPiquerezShaw25}
Edvard~Aksnes, Omid~Amini, Matthieu~Piquerez, and Kris~Shaw,
\emph{Cohomologically tropical varieties},
J. Inst. Math. Jussieu \textbf{24} (2025), no.~6, 2543--2572.
\url{https://doi.org/10.1017/S1474748025101114}.

\bibitem{AllermannRau10}
L.~Allermann and J.~Rau,
\emph{First steps in tropical intersection theory},
Math. Z. \textbf{264} (2010), no.~3, 633--670.

\bibitem{ArdilaKlivans06}
F.~Ardila and C.~Klivans,
\emph{The Bergman complex of a matroid and phylogenetic trees},
J. Combin. Theory Ser. B \textbf{96} (2006), no.~1, 38--49.

\bibitem{BabaeeHuh19}
F.~Babaee and J.~Huh,
\emph{A tropical approach to a generalized Hodge conjecture for positive currents},
Duke Math. J. \textbf{166} (2017), no.~14, 2749--2813.
\url{https://doi.org/10.1215/00127094-2017-0017}.

\bibitem{CoolsDraismaPayneRobeva12}
Filip~Cools, Jan~Draisma, Sam~Payne, and Elina~Robeva,
\emph{A tropical proof of the {Brill--Noether} theorem},
Adv. Math. \textbf{230} (2012), no.~2, 759--776.

\bibitem{Cox}
D.~A.~Cox, J.~B.~Little, and H.~K.~Schenck,
\emph{Toric Varieties},
Graduate Studies in Mathematics, vol.~124, Amer. Math. Soc., Providence, RI, 2011.

\bibitem{EsterovSCI}
Alexander~Esterov,
\emph{Sch{"o}n complete intersections},
preprint, 2024.
arXiv:,\href{https://arxiv.org/abs/2401.12090}{2401.12090} [math.AG].

\bibitem{EsterovECI}
Alexander~Esterov,
\emph{Engineered complete intersections: eliminating variables and understanding topology},
preprint, 2025.
arXiv:,\href{https://arxiv.org/abs/2504.16018}{2504.16018} [math.AG].

\bibitem{EsterovGusev}
A.~Esterov and G.~Gusev,
\emph{Systems of equations with a single solution},
J. Symbolic Comput. \textbf{68} (2015), no.~2, 116--130.

\bibitem{FeichtnerYuzvinsky04}
Eva~Maria~Feichtner and Sergey~Yuzvinsky,
\emph{Chow rings of toric varieties defined by atomic lattices},
Invent. Math. \textbf{155} (2004), no.~3, 515--536.
\url{https://doi.org/10.1007/s00222-003-0327-2}.

\bibitem{AlexFink}
A.~Fink,
\emph{Tropical cycles and Chow polytopes},
Beitr. Algebra Geom. \textbf{54} (2013), no.~1, 13--40.

\bibitem{FultonSturmfels97}
W.~Fulton and B.~Sturmfels,
\emph{Intersection theory on toric varieties},
Topology \textbf{36} (1997), no.~2, 335--353.

\bibitem{Katz12}
Eric~Katz,
\emph{Tropical intersection theory from toric varieties},
Collect. Math. \textbf{63} (2012), no.~1, 29--44.
\url{https://doi.org/10.1007/s13348-010-0014-8}.

\bibitem{Laz04}
R.~Lazarsfeld,
\emph{Positivity in Algebraic Geometry {I}: Classical Setting: Line Bundles and Linear Series},
Ergebnisse der Mathematik und ihrer Grenzgebiete. 3. Folge, vol.~48, Springer-Verlag, Berlin, 2004.

\bibitem{Li25}
Linxuan~Li,
\emph{Tropical fans supporting a reduced 0-dimensional complete intersection},
preprint, 2025, Version 2.
arXiv:,\href{https://arxiv.org/abs/2508.06694}{2508.06694} [math.CO].
\url{https://doi.org/10.48550/arXiv.2508.06694}.

\bibitem{MaclaganRincon18}
Diane~Maclagan and Felipe~Rinc{'o}n,
\emph{Tropical ideals},
Compos. Math. \textbf{154} (2018), no.~3, 640--670.
\url{https://doi.org/10.1112/S0010437X17008004}.

\bibitem{MaclaganSturmfels15}
D.~Maclagan and B.~Sturmfels,
\emph{Introduction to Tropical Geometry},
Graduate Studies in Mathematics, vol.~161, Amer. Math. Soc., Providence, RI, 2015.

\bibitem{Mikhalkin05}
Grigory~Mikhalkin,
\emph{Enumerative tropical algebraic geometry in {$\mathbb{R}^2$}},
J. Amer. Math. Soc. \textbf{18} (2005), no.~2, 313--377.
\url{https://doi.org/10.1090/S0894-0347-05-00477-7}.

\bibitem{MikhalkinRau}
G.~Mikhalkin and J.~Rau,
\emph{Tropical Geometry},
Oberwolfach Seminars, vol.~35, Birkh"auser, 2009.

\bibitem{MikhalkinZharkov14}
Grigory~Mikhalkin and Ilia~Zharkov,
\emph{Tropical eigenwave and intermediate {Jacobians}},
in \emph{Homological Mirror Symmetry and Tropical Geometry}, Lecture Notes of the Unione Matematica Italiana, vol.~15, Springer, Cham, 2014, pp.~309--349.
\url{https://doi.org/10.1007/978-3-319-06514-4_7}.
\bibitem{minkowski}
Hermann~Minkowski,
\emph{Theorie der konvexen K{"o}rper, insbesondere Begr{"u}ndung ihres Oberfl{"a}chenbegriffs},
in \emph{Gesammelte Abhandlungen}, vol.~2, Teubner, Leipzig, 1911, pp.~131--229.

\bibitem{Oxley}
James~Oxley,
\emph{Matroid Theory},
2nd ed., Oxford Graduate Texts in Mathematics, vol.~21, Oxford University Press, Oxford, 2011.

\end{thebibliography}
\end{document}